\documentclass{article}
\usepackage{amssymb}
\begin{document}

\title{Kahler geometry on toric manifolds, and some other manifolds with
large symmetry}
\author{S. K. Donaldson}
\maketitle


\tableofcontents

   \newtheorem{thm}{Theorem}
\newtheorem{lem}{Lemma}
\newtheorem{cor}{Corollary}
\newtheorem{prop}{Proposition}
\newtheorem{conj}{Conjecture}
\newcommand{\bC}{{\bf C}}
\newcommand{\bR}{{\bf R}}
\newcommand{\cH}{M}
\newcommand{\cP}{{\cal P}}
\newcommand{\bP}{{\bf P}}
\newcommand{\oa}{\overline{a}}
\newcommand{\ua}{\underline{a}}
\newcommand{\dbd}{\overline{\partial}}
\newcommand{\cA}{{\cal A}}
\newcommand{\SDiff}{{\rm SDiff}}
\newcommand{\cJ}{{\cal J}}
\newcommand{\cJint}{{\cal J}_{{\rm int}}}
\newcommand{\Riem}{{\rm Riem}}
\newcommand{\grad}{{\rm grad}}
\newcommand{\Ricci}{{\rm Ric}}
\newcommand{\ut}{\underline{t}}
\newcommand{\ux}{\underline{x}}
\newcommand{\utheta}{\underline{\theta}}
\newcommand{\bZ}{{\bf Z}}
\newcommand{\Xcx}{X_{{\rm cx.}}}
\newcommand{\Xsymp}{X_{{\rm symp}}}
\newcommand{\Xalg}{X_{{\rm alg}}}
\newcommand{\oP}{\overline{P}}
\newcommand{\Vol}{{\rm Vol}}
\newcommand{\Av}{{\rm Av}}
\newcommand{\cL}{{\cal L}}
\newcommand{\Lie}{{\rm Lie}}

\

\

\

\

In this article we discuss some aspects of existence problems in Kahler geometry;
a field which owes so much to  Yau. Kahler manifolds, in general,
are rather sophisticated mathematical objects and---the author feels---it
is often hard to acquire an intuition to back up the more abstract ideas.
Thus the main point  of this article is to discuss cases where
 the manifolds in question can be, to some extent, visualised and the existence
 problems stand out more clearly. Our central topic is the class of \lq \lq
 toric varieties'' and we will begin by reviewing the differential-geometric
 theory of these. Then we move on to consider  two variants of the toric
 condition---both involving manifolds with large symmetry groups---and make
 a special study of a Fano 3-fold found by Mukai and Umemura. In  the companion article (with R.S. Bunch), 
 immediately following this in the volume, we take  the \lq\lq visualisation''
 theme in a different direction, with numerical results for toric surfaces.

\

The author is grateful to Gabor Szekelyhidi, Rosa Sena-Dias and Yanir Rubinstein for pointing
out a substantial mistake in an earlier version of this paper. 
 
 \
 
\section{Background}
We begin by reviewing some basic notions in Kahler geometry. The author's
view of this subject is coloured by  an analogy with gauge theory
so, while it is only indirectly relevant, we will begin with that.

\subsection{Gauge theory and holomorphic bundles.}

Here we consider a complex vector bundle $E$ over a complex manifold $X$.
We want to study the interaction between two structures
\begin{itemize}
\item A hermitian metric on  $E$;
\item A holomorphic structure on $E$, which can be defined by
a $\dbd$-operator
$$   \dbd:\Omega^{0}(E)\rightarrow \Omega^{0,1}(E). $$
\end{itemize}

A basic fact is that given both of these structures there is a  unique compatible
unitary connection, in the sense that the $\dbd$-operator is the $(0,1)$-component
of the covariant derivative. 
Now there are two ways of setting up the theory. In the first---the traditional
point of view in complex geometry, as is \cite{kn:GH} for example---we fix
a holomorphic structure and  consider the various Hermitian  metrics. Then
we have, for example, the formula
\begin{equation}F_{h}= \dbd( h^{-1} \partial h)\end{equation}
for the curvature tensor in a local holomorphic trivialisation, where the
metric is defined by a  matrix-valued function $h$. In the second point of
view---closer
to what one does in general Yang-Mills theory---we fix the Hermitian metric
and consider various $\dbd$-operators.  We can identify the set of these
operators with the space
$\cA$ of unitary connections on $E$. This point of view brings in two infinite
dimensional groups. First, the group $U(E)$ of unitary automorphisms of $E$
 and second the group $GL(E)$ of general linear
automorphisms. Then $GL(E)$ acts on the space of $\dbd$-operators by conjugation,
and hence on the set $\cA$ of connections. The $\dbd$-operators which define
equivalent holomorphic structures are exactly those which are in the same
orbit of the $GL(E)$-action.

The advantage of this second point of view comes when studying the \lq\lq
jumping'' of holomorphic structures. This arises from the fact that the $GL(E)$
orbits are not usually closed in  $\cA$.  Fix a Kahler metric on the base space
$X$ and use this to define the Yang-Mills functional: the $L^{2}$ norm of
the curvature. When one seeks Yang-Mills connections compatible with a given
holomorphic structure ${\cal E}$ one attempts to minimise this functional over a $GL(E)$
orbit in $\cA$. But it may happen that there is no minimum, in the simplest
case because the infimum is achieved at a point in $\cA$ in the {\it closure}
but not in the orbit itself. Then one finds a Yang-Mills
connection not on the original holomorphic bundle ${\cal E}$, but on another
one ${\cal E}'$, such that there are arbitrarily small deformations of
${\cal E}'$ which are isomorphic to ${\cal E}$. This lies at the root
of  the link between Yang-Mills theory and  the {\it stability} of holomorphic bundles expressed by the Kobayashi-Hitchin conjecture \cite{kn:AB},
\cite{kn:UY},
\cite{kn:D1}.           

\subsection{Symplectic and complex structures}
Now we pass on to Kahler geometry. We study the interaction between
two structures on an underlying manifold $M$:
 a complex structure and  a symplectic form. We require these to be algebraically compatible in the sense
that the symplectic form is the imaginary part of a hermitian  metric. As before there
are two points of view we can take. In the first---the conventional point
of view in complex differential geometry---we fix the complex structure and 
vary the Kahler form. If we choose a reference form $\omega_{0}$ and vary
in the fixed cohomology class then (at least when $M$ is compact) any other
form can be represented by a Kahler potential, in the shape
$$  \omega_{\psi}= \omega_{0} + i \partial \dbd \psi. $$ 
For the alternative point of view we fix a symplectic form $\omega$ and consider
the space $\cJ$ of algebraically-compatible almost-complex structures on
$M$. Then the group $\SDiff$ of symplectomorphisms of $(M,\omega)$ acts on
$\cJ$, and this is the analogue of the  unitary gauge group $U(E)$ in the
previous case. We consider the subset $\cJint$ of integrable almost complex
structures, which is preserved by $\SDiff$. This is partitioned into equivalence
classes under the relation $J_{1}\sim J_{2}$ if $(M,J_{1}), (M,J_{2})$ are
isomorphic as complex manifolds. Although the group $\SDiff$ does not have
a true complexification one can argue that the equivalence classes in $\cJint$
are formally the orbits of such a (mythical) complexified group, in the sense
that they behave that way at the level of tangent spaces and Lie algebras
\cite{kn:D2}. 

\subsection{The equations}

The focus of this article is on the existence question for four different
kinds of special Kahler metrics, working within a fixed Kahler class on a
compact manifold.

\begin{enumerate}
\item {\it Extremal Kahler metrics}\ The definition is due to Calabi \cite{kn:Cal}. They are critical points (and in fact local minima) of the Calabi functional
$$  \int_{M} \vert \Riem(\omega) \vert^{2} d\mu_{\omega}, $$
where $\omega$ varies over the Kahler metrics in a fixed Kahler class and
$\Riem$
is the Riemann curvature tensor. The Euler-Lagrange equation is
$$  \dbd (\grad S_{\omega}) = 0,$$
where $\grad$ is the  gradient operator defined by $\omega$ and $S(\omega)$
is the scalar curvature. In other words,
the vector field $\grad S_{\omega}$ should be a holomorphic vector field.
On the face of it, this is a sixth order partial differential equation for
the Kahler potential $\psi$. 

\item {\it Constant scalar curvature Kahler metrics}\  These are just those with
$S_{\omega}$ constant.  Certainly they are extremal metrics (since the gradient
vanishes), and if it happens that $M$ has no non-trivial holomorphic vector
fields then an extremal metric must have constant scalar curvature.

\item {\it Kahler-Einstein metrics}\  By definition these are those where the Ricci
tensor  is  a multiple $\lambda \omega$. We will only consider the
case when $\lambda$ is positive (the zero and negative cases being completely
understood through the results of Yau and Aubin). By rescaling there is
no loss in supposing that $\lambda=1$. Solutions can only exist when $M$
is a \lq\lq Fano'' manifold and the class $[\omega]$ is $- c_{1}(M)$. 
 \item {\it Kahler-Ricci solitons}\ These again occur only in the Fano case.
 They are metrics for which
     $$  \Ricci - \omega= L_{v} \omega, $$
     where $L_{v}$ is the Lie derivative along  a holomorphic vector field $v$.
\end{enumerate}

Obviously a Kahler-Einstein metric has constant scalar curvature. There is
no simple relation between the other
two classes---extremal metrics and Kahler-Ricci solitons--- but they can
each  be thought of
as variants of the theory which take account of the possible holomorphic
vector fields on the manifold. All this is elucidated by the theory of the
{\it Futaki invariant}. We will not go in to this in detail here, since we
will see later how the theory works in explicit examples. Suffice it to say
that  in either situation the relevant holomorphic vector field which 
can be determined
{\it a priori} from standard topological data. More precisely, the vector
field it 
determined once we fix a maximal compact connected  subgroup of the 
group of holomorphic automorphisms. In either situation, an extremal metric
or Kahler-Ricci soliton will necessarily be Einstein/constant scalar curvature
if the Futaki invariant vanishes.

There is, of course, as yet no general existence theory for these structures
but at the conjectural level one can see a detailed analogy with the Yang-Mills
case. We do not want to go into this further here---partly because
there is a comprehensive recent survey article \cite{kn:PhongSturm}---but proceed with our study
of special classes of manifolds.

\section{Toric manifolds}

We say  that a compact Kahler manifold $X$ of complex dimension $n$ is {\it toric} if the compact torus $T^{n}$ acts by isometries on $X$ and the
extension of the action to the complex torus $T^{n}_{c} \cong ( \bC^{*})^{n}$
acts holomorphically with a free, open, dense orbit $X_{0}\subset X$.

\

\subsection{Local differential geometry}
\subsubsection{Complex coordinates}
Here we work in the neighbourhood of a point in the free orbit $X_{0}$. We
can use the group action to define local co-ordinates. So  we have complex
co-ordinates
$$  \tau_{a} =\frac{1}{2}( t_{a} + i \theta_{a})$$
say. The factor $2$ here will simplify the formulae later. Locally the isometry group acts by translations
in the $\theta_{a}$ directions. (Later, when we work globally, the $\theta_{a}$
will become \lq\lq angular'' co-ordinates, with period $4\pi$.) Locally,
a Kahler metric is given by $i\partial \dbd \phi$ for a function $\phi$ of
the complex variables $\tau_{a}$. If this function only depends on the real
parts $t_{a}$ then the metric will obviously be invariant under translations
in the $\theta_{a}$ directions and it is not hard to see that any metric
of the kind we are considering arises in this way. Now if we write
$\phi=\phi(t_{a})$ then the tensor $i\partial \dbd \phi$ is just
$$   \sum_{a b} \frac{\partial \phi}{\partial t_{a} \partial t_{b}} d\tau_{a}
d\overline{\tau}_{b}, $$
and this defines a positive Hermitian form if and only if the Hessian matrix
of $\phi$ is positive definite; or in other words $\phi$ is a {\it convex}
function of the real variables $\tau_{a}$. Thus the theory of
convex functions on Euclidean spaces is embedded, as this translationally
invariant case, in the theory of Kahler geometry. We write $\nabla^{2} \phi$
for the Hessian of $\phi$ and also use index notation $\nabla^{2}\phi=(\phi^{ab})$.
The placing of the indices is unconventional but will be convenient later.
We write
$(\phi_{ab})$ for the inverse matrix.
 Explicitly the symplectic
form $\omega$ is
  $$  \frac{1}{2}\sum \phi^{ab} dt_{a}\wedge d\theta_{b},
  $$ and the Riemannian metric is
  $$ \frac{1}{2} \left( \sum \phi^{ab} dt^{a} dt^{b} + \sum \phi^{ab} d\theta^{a} d\theta^{b}\right).
  $$
  We regard the curvature tensor of this metric as an element of $\Lambda^{2}
  \otimes \Lambda^{2}$. Then the curvature tensor is
  $$   \sum R^{abcd} d\tau_{a} d\overline{\tau}_{b} \otimes d\tau_{c}d\overline{\tau}_{d},
  $$
  where
  \begin{equation} R^{abcd} = \phi^{abcd} - \phi^{ac \lambda} \phi^{bd\mu} \phi_{\lambda
  \mu}. \end{equation}
  
  (Here we use the summation convention over the repeated indices.
The third and fourth order derivatives of $\phi$ are written as $\phi^{abc},
\phi^{abcd}$ in the obvious way.)  This formula for the curvature
tensor is just the formula (1), expressed in our current notation.)

\subsubsection{Symplectic coordinates}
We now take a different point of view, following Guillemin \cite{kn:Guil} and
Abreu \cite{kn:Ab}, and also the general scheme outlined in the previous section.
Thus we consider an open set in $\bR^{n}\times \bR^{n}$ with linear coordinates
$x^{a}, \theta_{a}$. More invariantly, we should write the ambient space
as $V\times V^{*}$ where $V=\bR^{n}$, with coordinates $x^{a}$. We assume
the open set has the form $Q\times V^{*}$ where $Q\subset V$ is convex. On this open
set we consider the standard symplectic form
$$\Omega= \frac{1}{2} \sum dx^{a} d\theta_{a}. $$
This is preserved by the translations in the $\theta$ variables. More precisely
we have a Hamiltonian action of the group $G=V^{*}$ on the symplectic manifold
$Q\times V^{*}$ and the moment map is just the projection to $Q$, with components
the coordinates $x^{a}$. We consider $G$-invariant almost-complex structures
on $Q\times V^{*}$, algebraically compatible with $\Omega$. Now at each point
such a structure is specified by a subspace of the complexified cotangent
bundle which has a unique basis of the form
$$  \epsilon_{a} = d\theta_{a} + Z_{ab} dx^{b}, $$
where $(Z_{ab})$ is a symmetric complex matrix with positive definite imaginary
part. (This is just the standard description of the Siegel upper half-space
$Sp(n,\bR)/U(n)$.) So our almost-complex structure is represented by a matrix-valued
function $(Z_{ab})$ and $G$-invariance specifies that $Z$ is a function of
the variables $x^{a}$. Following our general scheme we should now determine
when such an almost-complex structure is integrable. By definition this means
that the $2$-forms 
$$d\epsilon_{a}= \frac{\partial Z_{ab}}{\partial x^{c}} dx^{c} dx^{b} $$ can be expressed as $\sum \alpha_{a b}
\wedge \epsilon_{b}$ and  this only happens when all the 
$d\epsilon_{a}$ are zero (since $d\epsilon_{a}$ does not contain any terms
involving $d\theta_{i}$). So the integrability condition is 
\begin{equation}  \frac{\partial Z_{ab}}{\partial x^{c}}=\frac{\partial Z_{ac}}{\partial
x^{b}}. \end{equation}
Now consider the action of the infinite-dimensional symplectomorphism group.
In this situation we need to consider the symplectic diffeomorphisms that
commute with the $G$-action. More precisely we want to take the Hamiltonian
diffeomorphisms generated by functions that Poisson-commute with the generators
of the $G$-action; but  these are just the functions of the $x^{i}$ variables.
The corresponding group ${\cal G}$ of diffeomorphisms can be identified with smooth functions on $Q$,  where a function $f$  acts by taking a point $(\ux,\utheta)$ to
$(\ux, \utheta+ D f)$. This gives an action on the space of almost-complex
structures which simply takes $Z_{ab}$ to $Z_{ab} + f_{ab}$, where $f_{ab}$ is the Hessian
of $f$. 

Now consider  the action of ${\cal G}$ on the integrable structures. The
condition (3) implies, by the elementary \lq\lq criterion for an exact differential'', that there are complex-valued functions $i t_{a}$
such that
$$   Z_{ab}= i \frac{\partial t_{a}}{\partial x^{b}}. $$

The fact that $Z_{ab}$ is symmetric implies, by the same criterion, that
there is a single complex valued function $F$ such that $t_{a}=\frac{\partial
F}{\partial x^{a}}$, in other words
$$  Z_{ab}= \frac{\partial^{2} F}{\partial x^{a}\partial x^{b}}. $$
If we let $f$ be minus the real part of $F$ then the action of $f\in
{\cal G}$ takes the structure $(Z_{ab})$ to a new structure with zero real part.
So ,taking account of this diffeomorphism group, we can reduce to considering $Z=iY$, with $Y$ real and positive definite. Now the functions $t_{a}$
are real and $\epsilon_{a}= d( t_{a}+ i \theta_{a})$ so
$t_{a}+ i\theta_{a}$ are local complex co-ordinates. (Thus we confirm
the Newlander-Nirenberg integrability theorem in this special case.).  Write $u$ for the
imaginary part of the function $F$ above, so $$Y_{ab}= \frac{\partial^{2}
u}{\partial x^{a} \partial x^{b}}= u_{ab} . $$

Some linear algebra shows that the metric defined by the almost complex structure
and the fixed form $\Omega$ is
\begin{equation}  \frac{1}{2} \sum u_{ij} dx^{i} dx^{j} + u^{ij} d\theta_{i} d\theta_{j}, \end{equation}
where $(u^{ij})$ is the matrix inverse of the Hessian $(u_{ij})$.

\

The conclusion of this is that we have another description of the local differential
geometry, defined by a convex function
$u$ of the variables $x^{a}$. The relation between this picture and that
in complex co-ordinates discussed above is just the {\it Legendre transform}
for convex functions. That is, given a convex function $u$ on $Q\subset V$
we define a function $\phi$ on an open set $Q^{*}\subset V^{*}$ by decreeing
that 
$$ \phi(\ut)= \sum x^{a}t_{a}- u(\underline{x}), $$
where the point $\underline{x}\in V$ is the unique point where $Du=\ut$.
As is well-known, this transform expresses a symmetric relation between $u$
and $\phi$, so $u$ is the Legendre transform of $\phi$. Further, the Hessian
$\phi^{ab}=\frac{\partial^{2} \phi}{\partial t_{a}t_{b}}$ is the inverse of the
Hessian $u_{ab}$ of $u$ at the corresponding point. It is easy to see using
this that the Legendre transform does give a Kahler potential for the same
metric expressed in the complex co-ordinates. Conversely if we start with
the complex description and a convex function $\phi$ then  the Legendre transform
gives the symplectic picture. More invariantly, the map $\ux$ is characterised
as the moment
map for the action of the group of translations. 

Thus we have two natural coordinate systems to use when discussing this local
differential geometry, and of course we can transform any formulae from one
set-up to the other. Working in the symplectic picture we set
$$  F_{ij k l} = u_{ia} u_{jb} \frac{\partial^{2} u^{ab}} {\partial x^{k} \partial
x^{l}}.$$
Then one finds that the Riemann curvature tensor is
\begin{equation}  F_{ijkl} \eta^{i}\wedge \eta^{k}\otimes \eta{j}\wedge \eta^{l},\end{equation}
where $\eta^{a}=dx^{a}+ i u^{ab} d\theta_{b}$.
So the four-index tensor $F$ is essentially the same as the curvature tensor.
For example the norm if the Riemann curvature tensor is the same as the natural
norm of $F$ i.e.
$$ \vert F \vert^{2} = \sum F_{ijkl} F_{abcd} u^{ia} u^{jb} u^{kc} u^{l d}.
$$
The Ricci tensor is in the same fashion, equivalent to the tensor
$$   G_{ij}= F_{ijkl}u^{kl} , $$
which can also be expressed as $$G_{ij}=\frac{\partial^{2} L}{\partial x^{i} \partial
x^{j}}$$
where $L=\log \det(u_{ij})$. The scalar curvature is given by another contraction
yielding Abreu's   formula
\begin{equation} S=  G_{ij}u^{ij}= \sum_{ij} \frac{\partial^{2} u^{ij}}{\partial x^{i}\partial
x^{j}}. \end{equation}

\

We mentioned in the previous section that in the general case the group of
symplectomorphisms does not have a complexification, and this limits the
practicality of the symplectic approach to  Kahler geometry. But in this special
situation there {\it is} a complexification of ${\cal G}$:  simply the complex
valued functions on $Q$ under addition. Further, in it is nearly true that this complexified group ${\cal
G}^{c}$  acts on the set of almost complex structures, represented as matrix-valued
functions $Z_{ab}$. The \lq\lq action'' is simply to map $Z$ to $Z+\frac{\partial^{2}
F}{\partial x^{a}\partial x^{b}}$. It is only a local action because
the condition that the imaginary part of $X$ is positive definite could be
violated. Our discussion above  asserts that
all the integrable structures are in a single orbit of this complexified
action and the parametrisation by the function $u$ is  the parametrisation
by an open set in the quotient ${\cal G}^{c}/{\cal G}$.  Further, it is easy to verify in this framework
that the scalar curvature given by the formula (6) is a moment map for the
action of ${\cal G}$ with respect to the natural symplectic structure on the space
of almost-complex structures (which is derived from the invariant symplectic
form on the Siegel upper half space), see \cite{kn:D4}.

\subsection{The global structure}

In the previous section we discussed the local differential geometry of a
toric manifold in the dense open set where the torus action is free. We now
go on to the
global picture. There are at
least three different points of view we can take but the essential thing is that
this structure is encoded by a bounded polytope $P\subset \bR^{n}$, or more invariantly
$P\subset V$ in the notation of the previous section. This polytope is defined
by a finite collection of linear inequalities $\lambda_{r}(\ux)>c_{r}$ corresponding
to the codimension-$1$ faces. So $\lambda_{r}$ are vectors in the dual space
$V^{*}$. We suppose that there is an integer lattice in $V$, which we can
take to be the standard $\bZ^{n}$ in $\bR^{n}$. Then there is a dual lattice
in $V^{*}$ and we suppose that the $\lambda_{r}$ lie in this dual lattice.
We can rescale so that the $\lambda_{r}$ are primitive vectors with respect
to this lattice. 
Further, we suppose that  each vertex of $P$ is contained in exactly $n$
codimension faces and that the corresponding $\lambda_{r}$ form an integer
basis for the dual lattice.  Such
a polytope  is called a {\it Delzant polytope}. 
Another way of expressing the condition is via the group $\Gamma$ of maps
$$ \ux\mapsto A \ux + \underline{b} $$
from $\bR^{n}$ to itself, where $A$ is restricted to lie in $GL(n,\bZ)$. Up
to the action of $\Gamma$, a neighbourhood of any vertex of $P$ is equivalent to a neighbourhood
of $0$ in the infinite polytope $\{x_{i}>0\}\subset \bR^{n}$. If the vertices of the polytope
are integral we call it an integral Delzant polytope.

\

 {\bf Example}\  The standard simplex in $\bR^{n}$, given
by the inequalities
$$x^{1}>0, x^{2}>0, \dots, x^{n}>0, x^{1}+x^{2} + \dots x^{n}\leq 1$$  
is a Delzant polytope.
\subsubsection{Complex charts}

Start with a Delzant polytope $P$. Let ${\cal S}$ be the finite set
of pairs of
\begin{itemize}
\item a vertex $p$ of $P$;
\item an ordering $\lambda_{r(i)}$ of the faces containing $p$.
\end{itemize}
For any two $\sigma=(p,r(\ ))$ and $\sigma'=(p', r'(\ ))$ in ${\cal S}$ there is a unique element $\gamma_{\sigma, \sigma'}$ of $\Gamma$ which maps $p$
to $p'$ and matches up the corresponding faces. Obviously we have

 $$\gamma_{\sigma, \sigma}= 1\ ;\ 
\gamma_{\sigma, \sigma'}= \gamma_{\sigma',\sigma}^{-1}\ ;\ 
\gamma_{\sigma, \sigma''}= \gamma_{\sigma, \sigma'}\circ \gamma_{\sigma'
\sigma''}. $$
Now suppose we have any space $M^{*}$ on which $\Gamma$ acts and $M^{*}$ is a
subset of a larger space $M$.We take the product
${\cal S}\times M$ and define a relation 
$$  (\sigma, m)\sim (\sigma', \gamma_{\sigma,\sigma'} (m)), $$
for $m\in M^{*}$.
The properties above tell us that this is an equivalence relation, so we can
take the quotient ${\cal S}\times M/\sim$.
In our case we take $M$ to be $\bC^{n}$ and $M^{*}=(\bC^{*})^{n}\subset \bC^{n}$. Then $GL(n,\bZ)$ acts on $M^{*}$.
This is clear if we identify $\bC^{*}$ with $\bC/\bZ$ and hence $M^{*}$ with
$\bC^{n}/\bZ^{n}$. In terms of the original description, with co-ordinates
$z_{i}$ on $\bC^{n}$, we make a matrix $(a_{ij})$
act on $(\bC^{*})^{n}$by
$$ z'_{i} = \prod z_{j}^{a_{ij}}, $$
which is well-defined since the $a_{ij}$ are integers. There is a natural
homomorphism from $\Gamma$ to $GL(n, \bZ)$ so $\Gamma$ acts on $M^{*}$ via
this. Then it is clear from
the construction that  the quotient $\Xcx$ is a complex manifold covered by charts $M_{\sigma}$
labelled
by elements of $\Sigma$, each chart being a copy of $M=\bC^{n}$. The charts
for the $n!$ different elements of $\Sigma$ belonging to the same vertex
of $P$ have the same image so it suffices just to take one of them.  There
is an   action of the
complex torus $T^{n}_{c}$ with a dense orbit, which is the image of any
$\{\sigma\}\times M^{*}$. The construction behaves well with respect to restriction
to faces, so for each $m$-dimensional face $\Pi$ of $P$ there is a submanifold
$X^{\Pi}\subset \Xcx$ which is an $m$-dimensional complex submanifold with an
action of $T^{m}_{c}$ induced from the action on $\Xcx$. Indeed the orbits of
the $T^{n}_{c}$ action on $\Xcx$ correspond to these faces. In particular the vertices
of $P$ correspond to points of $\Xcx$; the fixed points under the $T^{n}_{c}$
action.

\

{\bf Example} When $P$ is the $n$-simplex, as  above, the manifold $\Xcx$ we construct
is $\bC\bP^{n}$.

\

So far we have not used the full strength of the data we began with. For
example, we could simply have omitted some vertices of $P$ and run the same
construction. We have also  thrown away some of the data, through the
homomomorphism from $\Gamma$ to $GL(n,\bZ)$.  First, the fact that the vertices come from a bounded polytope yields
the compactness of the space $\Xcx$ we have defined. We leave this as an exercise
for the reader.  For the second point, it is indeed the case that if we vary
the constants $c_{r}$ slightly (so that we do not introduce or remove any
vertices) we get the same complex manifold $\Xcx$. The extra structure of the
specific polytope corresponds
to fixing a distinguished cohomology class in $H^{2}(\Xcx;\bR)$. This is easiest to see in
the case when the polytope is integral. Then the $\gamma_{\sigma \sigma'}$
lie in a smaller group $\Gamma_{\bZ}\subset \Gamma$ which is an extension
$$   \bZ^{n}\rightarrow \Gamma_{\bZ} \rightarrow GL(n,\bZ). $$
We take the trivial complex line bundle $\underline{\bC}$ over $M=\bC^{n}$.
Then $\Gamma_{\bZ}$ acts on the restriction of $\underline{\bC}$ to $M^{*}$
and the same construction gives a complex line bundle $L\rightarrow \Xcx$. Furthermore
this is an equivariant line bundle for the  $T^{n}_{c}$ action. The distinguished
cohomology
class is just the first Chern class of $L$. In general, when the vertices
are not integral we consider the sheaf $Z^{1}$ of closed $1$-forms over $\Xcx$. We can use the $\gamma_{\sigma \sigma'}$ to define a closed $1$-form on $M_{\sigma}\cap
M_{\sigma'}$ and this yields a Cech cocycle with values in this sheaf. Then
the short exact sequence of sheaves
$$    0 \rightarrow \bR \rightarrow C^{\infty}(\Xcx) \rightarrow Z^{1} \rightarrow
0$$
gives a boundary map from $H^{1}(\Xcx; Z^{1})$ to $H^{2}(\Xcx,\bR)$ which defines
the distinguished cohomology class. (In fact this cohomology class is not
changed if we translate $P$. A more precise statement is that the Delzant polytope
$P$ can be recovered from the complex manifold $X$ with a suitable distinguished $T^{n}_{c}$-{\it
equivariant} cohomology class.) 

\

{\bf Example.} Consider a vertex $p$ of a Delzant polytope $P$. There is no
loss of generality in supposing that $p$ is the origin and that near the
origin $P$ agrees with the standard model $\{x^{i}>0\}$. Then, for $\delta>0$,
we define $P_{\delta}$ to be the subset of $P$ defined by the additional
inequality $\sum x_{i}> \delta$.  For small enough $\delta$ this is again
a Delzant polytope and the complex manifold $X_{\delta}$ is the {\it blow-up}
of $X$ at the fixed point corresponding to $P$. The exceptional divisor $E$
is a copy of projective space, associated to the \lq\lq new'' $n-1$-simplex in the boundary of $P_{\delta}$.  The manifold does not
vary with $\delta$ but the evaluation of the distinguished cohomology class on the standard generator of $H_{2}(E)\subset H_{2}(X_{\delta})$ is $\delta$.

\

Now we  go back to differential geometry. If we have a Kahler metric on
$\Xcx$, its restriction to the open orbit is described by a Kahler potential;
a convex function $\phi$ on $\bR^{n}$, as above. Conversely we can define
am \lq\lq admissible'' convex function $\phi$ to be one which defines a Kahler
metric over the orbit which extends smoothly to the compact manifold. This
is a condition on the asymptotic behaviour of $\phi$ at infinity in $\bR^{n}$.
The essence of the condition is that $\phi$ is asymptotic to the piecewise
linear function
$$   \Phi(\ut) = \max_{p} p.\ut, $$
where $p$ runs over the vertices of the polytope. Thus if we let
$\phi_{\lambda}$ be the rescaling $\phi_{\lambda}(\ut)= \lambda^{-1} \phi(\lambda
\ut)$ for $\lambda \in \bR$ then $\phi_{\lambda}\rightarrow \Phi$ (in $C^{0}$)as $\lambda$
tends to infinity. In the model case when $0$ is a vertex and $P$ agrees
locally with $\{ x^{i}>0\}$ the local complex co-ordinates are $z_{a}= \log
\tau_{a}$ and so $\vert z_{a}\vert^{2}= e^{t_{a}}$. The admissible condition
is that $\phi$ extends to a smooth function of the complex co-ordinates $z_{a}$.

\

{\bf Example}
The round metric on the $2$-sphere with area $2\pi$ is given by the Kahler potential
$$  \phi(t)=  \log (1+e^{t}). $$
In terms of a local complex co-ordinate $z$ this is $ \log (1+
\vert z\vert^{2})$.

\subsubsection{Symplectic construction}

Here we start with the product $P\times T^{n}$ with  standard co-ordinates
$x^{a}, \theta_{a}$ as before, except of course that now the $\theta_{a}$
are taken to be \lq\lq angular'' co-ordinates with period $4\pi$. This is a noncompact
symplectic manifold with the standard symplectic
form $\Omega= \sum dx^{a} d\theta_{a}$  and with Hamiltionian $T^{n}$ action whose moment map is the projection to  $P$. The essential point is that this can be compactified to a compact
symplectic manifold $\Xsymp$ and the moment map extends to a map with image
the closure $\overline{P}$. This works in a similar fashion to the complex
picture. For example, consider the neighbourhood of a vertex of $P$ which
as usual we can take to be the origin, with $P$ locally modelled on $\{x^{i}>0\}$.
Then $\Omega$ is the pull-back of the standard form on $\bC^{n}$ under the
map
$$  (x^{a}, \theta_{a})\mapsto  (\vert x_{a} \vert^{1/2} e^{i\theta_{a}}),
$$
 We adjoin a neighbourhood
of $0$ in $\bC^{n}$ to $P\times T^{n}$ using this map and repeat the construction,
modified in the obvious way, for all other boundary points of $P$. 

Now of course this symplectic construction describes the same object as
the complex construction in the previous section. We return to the discussion
of the local differential geometry taking now $Q=P$. We can start with an
admissible Kahler potential  $\phi$ on $\bR^{n}= V^{*}$. Then its Legendre transform
is a function on $P$. Around a vertex, as above,  this has the form
$$  u= \sum x^{i} \log x^{i} + v , $$
where $v$ is a smooth function (on the manifold with corners). We say that
a symplectic potential $u$ is admissible if it is the Legendre transform
of an admissible Kahler potential $\phi$. Stated explicitly in terms of $u$
this the requirement of
\lq\lq Guillemin boundary conditions'', which are
\begin{enumerate}
\item $u$ is a continuous function on $\oP$, smooth in the interior.
\item The restriction of $u$ to each face is smooth and strictly convex.
\item Let $q$ a boundary point which lies on a codimension $r$ face of $P$,
so without loss of generality $q =0$ and $P$ is locally defined by equations
$x^{1}>0, \dots x^{r}>0$. Then near $q$
$$   u= \sum_{i=1}^{r} x_{i} \log x_{i} + v $$
where $v$ is smooth.
\end{enumerate}

It is easy to see that such functions exist. For example we can take the
Guillemin function
$$  u= \sum_{r} (\lambda_{r}- c_{r}) \log (\lambda_{r}-c_{r}). $$
Either way, we get a map from the complex manifold $\Xcx$ to the symplectic manifold
$\Xsymp$ which matches up the structures involved.

{\bf Example} The round metric on $S^{2}$, of area $2\pi$,  is defined by the symplectic
potential, on the interval $[0,1]$,
$$  u(x)=\left( x \log x + (1-x) \log (1-x)\right). $$

\subsubsection{Algebraic construction}
Here we suppose that the Delzant polytope $P$ is integral. We consider all
the multiples $k\oP$ for integers $k\geq 0$ and let $B_{k}$ be the set of lattice
points
$$  B_{k}= k \oP \cap \bZ^{n}. $$
Let the number of points in $B_{k}$ be $N_{k}+1$. We can put all these sets
together by considering the cone over $P$
$$ cone(P)= \{(\ux, y)\in \bR^{n+1}: y\geq 0,  \ux\in y \oP\}. $$
 The disjoint union of the sets $B_{k}$ can be identified with the set $B=cone(P)\cap \bZ^{n+1}$.
Now $B$ is an abelian semi-group under addition and we have a corresponding
ring $R$ over $\bC$ with one generator $s_{b}$ for each point of $b\in B$
and relations
   $s_{b} s_{b'}= s_{b+b'}$.
   This is a graded ring, $R=\bigoplus R_{k}$, where
   $R_{k}$ has a basis $s_{\nu}$ corresponding to the points $\nu$ of $B_{k}$.   Further, there is an obvious action of the torus $T_{c}^{n}$ on $R$.

   All of these definitions make sense for any convex set $P$. The crucial
   fact is that when the $P$ is an integral polytope the ring is {\it finitely
   generated}. Thus there is a corresponding projective variety $\Xalg={\rm
   Proj}(R)$,
   and the group action on $R$ defines an action on $\Xalg$. Second, if $P$
   is Delzant, then $\Xalg$ is smooth and of course this recovers the same
   complex manifold $\Xcx$. The vector spaces $R_{k}$ are the  sections
   $$ R_{k}= H^{0}(\Xcx, L^{k})$$
   and it is not hard to see  that for any $k\geq 1$ the sections give an
   embedding $\Xcx\rightarrow \bP(R_{k}^{*})$. From this algebro-geometric
   point of view the integer $\lambda_{r}(\nu)-c_{r}$, for lattice points $\nu\in
   \oP$, is the order of vanishing of the section $s_{\nu}$ along the corresponding
   divisor
   in $\Xcx$.
   
   \
   
   {\bf Example}
   Let $P$ be the square $(0,1)^{2}\subset \bR^{2}$. The corresponding manifold
   is the product $S^{2}\times S^{2}$. The points in $B_{1}$ are the four
   vertices $p_{0}=(0,0), p_{1}=(0,1), p_{2}=(1,0), p_{3}=(1,1)$ so $R_{1}$ has a corresponding basis 
   $s_{0} s_{1}, s_{2}, s_{3}$ say. The equation $p_{0}+p_{3}= p_{1}+p_{2}$
   goes over to the relation $s_{0} s_{3}= s_{1}s_{2}$. The embedding of
   $\Xcx$ in $\bP^{3}$ has image the quadric hypersurface cut out by the
   equation $Z_{0} Z_{1}- Z_{2} Z_{3}=0$.
   
   \

   When the polytope $P$ is integral but not Delzant the variety $\Xalg$
   we construct is singular. If each vertex lies on exactly  $n$ codimension-1
    faces then $\Xalg$ is an orbifold. Much of the theory, including the differential-geometric constructions, extends easily to this case.

   \
   
   To sum up we have three ways---complex, symplectic and algebraic--- of constructing a compact manifold associated to an integral Delzant polytope.
From now on we will just denote this by $X$.

\

   \subsubsection{Real forms}
   A toric manifold $X$ contains a  submanifold $X_{\bR}$ of one half the dimension which is
   a \lq\lq real form'' in the complex picture and Lagrangian in the symplectic
   picture. To define this from the first point of view we just observe that
   the action of $\Gamma$ on  $M^{*}= (\bC^{*})^{n}$ preserves the subset
   $M_{\bR}^{*}$ of real points. Then we run the same construction. From
   the symplectic point of view we let $A$ be the subgroup of the real torus
   $T^{n}$  given by the elements of order $2$, so $A$ is isomorphic to $(\bZ/2)^{n}$.
   Then we consider the subset $A\times P\subset T^{n}\times P$ and check
   that the closure of this in $X$ is a smooth $n$-dimensional manifold.
   From the algebro-geometric point of view we simply observe that all our
   relations are real, so complex conjugation acts on everything and we get
   a real form of our complex algebraic variety.
   
   This construction is particularly vivid in the symplectic picture \cite{kn:Guil2}.
   The
   composite
   $$   X_{\bR}\rightarrow X \rightarrow \oP, $$
   is a $2^{n}$-fold covering map over the interior $P\subset \oP$ so we
   can construct $X_{\bR}$ by taking $2^{n}$ copies of $\oP$ and gluing the
   boundary components appropriately. The Riemannian metric on $P$ given
   by the Hessian $u_{ij}$ of an admissible symplectic potential extends to a smooth
   Riemannian metric on $X_{\bR}$. In particular we get a conformal
   structure on $X_{\bR}$ and when $n=2$ a Riemann surface structure on the
   oriented cover of $X_{\bR}$. (The surface $X_{\bR}$ is only itself orientable
   in the case when $P$ is a rectangle.) For example, if $P$ is the standard
   triangle in $\bR^{2}$ then $X_{\bR}$ is a real projective plane in
   $X=\bC\bP^{2}$ and can be constructed by gluing four triangles. The
   oriented cover is $S^{2}$, constructed by gluing eight triangles. In general
   we get a class of Riemann surfaces obtained by gluing eight polygons.
   Given a symplectic potential $u$, the induced conformal structure on $\oP$ is equivalent
   to the standard disc. So if $P$ has $s$ vertices we get an invariant of
   $u$ in the moduli space ${\cal M}_{s}$ of configurations of $s$ distinct points on $S^{1}=\bR\bP^{1}$
   modulo the action of $PSL(2,\bR)$. This determines the conformal structure
   of $X_{\bR}$, and is an interesting global invariant of a toric Kahler
   surface.

   \subsection{Algebraic metrics and asymptotics}

  If $X$ is any compact complex manifold and $L\rightarrow X$ a very ample line
  bundle we can generate Kahler metrics on $X$ by the following procedure.
  Choose a Hermitian metric on the complex vector space $H^{0}(X;L)$. This
  induces a metric on the dual space and hence a standard Fubini-Study metric
  on the complex projective space $\bP(H^{*}(X;L)^{*})$. Now we use the embedding
  $\iota:X\rightarrow \bP(H^{0}(X;L)^{*})$  to induce a Kahler metric on
  $X$. We  call metrics of this kind \lq\lq algebraic Kahler metrics''.
  
  This construction becomes very simple and explicit in the toric case. We
  consider metrics on $H^{0}(L)$ which are invariant under the torus action,
  hence are diagonal in the standard basis $s_{\nu}$.  A collection
  of positive numbers
  $a_{\nu}$, for each lattice point $\nu$ in $\oP$, defines an invariant metric
  with $\Vert s_{\nu}\Vert^{2}= a_{\nu}^{-1}$. Given this data $\{a_{\nu}\}$
  we have
  a Kahler potential on $\bR^{n}$:
  \begin{equation} \phi(\ut)= \log \left( \sum_{\nu} a_{\nu} e^{\nu. \ut}\right), \end{equation}
  where $\nu.t$ denotes the dual pairing between the copy of $\bR^{n}$ on
  which $\phi$ defined and the copy of $\bR^{n}$ containing $P$. This is the potential
  which defines the algebraic metric via the projective embedding.

  We will not discuss this topic at length here, but we want to make the
  point
  that the data $ -\log a_{\nu}$---a real-valued function on the lattice
  points in $\oP$---can be thought of as a \lq\lq discrete approximation''
  to the {\it symplectic potential} $u$---a real-valued function on $\oP$.
  This only makes sense as an asymptotic statement, when we replace the bundle
  $L$ by $L^{k}$ and $P$ by $kP$ for large $k$. Rescaling, we can equivalently fix $P$
  and replace the integer lattice by $k^{-1} \bZ^{n}$. We discuss two simple
  precise statements which illustrate this general idea but for many further developments in a similar vein
  we refer to the recent works of Zelditch \cite{kn:Zeld}.
 
 \

  \subsubsection{ Asymptotics of $L^{2}$-metrics}

  Suppose we start with some symplectic potential $u$ and corresponding Kahler
  potential $\phi$. Then $\phi$ can be regarded as a Hermitian metric on
  the line bundle $L$ over the toric variety. Thus we have a natural $L^{2}$-metric
  on $H^{0}(X;L)$
  $$  \Vert s \Vert^{2} = \int_{X} \vert s\vert^{2} d\mu_{\phi}, $$
  where the pointwise norm $\vert s \vert$ is defined by $\phi$ and $d\mu_{\phi}$
  is the volume form of the Kahler metric. Thus, starting with $u$ we get
  a collection of numbers $a_{\nu}= \Vert s_{\nu} \Vert^{-1}$. Now replace $L$ by $L^{k}$, as above.
  The same symplectic potential $u$ defines a metric on $L^{k}$ and we get
  a collection of numbers $a_{\nu}^{(k)}$ say, for $\nu\in \oP\cap k^{-1}
  \bZ^{n}$. 
  One precise statement expressing the general idea above is that for
  each $\epsilon>0$ and compact subset $K\subset P$ there is a $k_{0}$
  such that
  $$     \vert u(\nu)-k^{-1} \log a_{\nu}^{(k)}\vert <\epsilon, $$
  once $k\geq k_{0}$, for all $\nu\in K\cap k^{-1} \bZ^{n}$.
  
  The proof of this is very simple. Go back to the case $k=1$ for the moment.
  Unravelling the definitions, the coefficients $a_{\nu}$ are given by
  $$  a_{\nu}^{-1} = \int_{\bR^{n}} e^{-\phi} e^{ \ut.\nu}   \det(\nabla^{2}
  \phi)\ 
  d\ut, $$
  where $\phi$ is the given Kahler potential. (Notice, by the way, that Holder's inequality shows that $\nu \mapsto -\log a_{\nu}$ is a convex function, in
the obvious sense.)
  Rescaling, we get
  $  a_{\nu,k}^{-1}= I_{\nu}(k)$ say, where
  \begin{equation} I_{\nu}(k)= \int_{\bR^{n}} e^{-k ( \phi- \ut.\nu)} \det(\nabla^{2}\phi)
  d\ut. \end{equation}
  (Notice that these formulae make sense for {\it any} $\nu\in \oP$ and the restriction
  to the lattice $k^{-1} \bZ^{n}$ is not really relevant here.) So we see
  that our question reduces to  the standard discussion of the asymptotic
  behaviour of the integral * as $k\rightarrow \infty$. The dominant contribution comes from the a neighbourhood of the point $\ut_{0}$ where $\phi-\ut.\nu$ is
minimal and the standard Laplace approximation is
$$  I_{\nu}(k) \sim (2\pi k)^{-n/2} {\rm exp}( -k(\phi(t_{0})-t_{0} \nu))
\det\nabla^{2}\phi (\ut_{0}). $$
But $\ut_{0}$ is just the point which corresponds to $\nu$ under the Legendre
transform, and $\phi(\ut_{0})-\ut_{0}.\nu$ is $-u(\nu)$. So
$$  k^{-1} \log I_{\nu}(k) = u(\nu) + O(k^{-1} \log k), $$
and our result follows since $k^{-1}\log k\rightarrow 0$ as $k\rightarrow
\infty$.

\

Following on this line, it is easy to derive a special case of Tian's Theorem
from \cite{kn:T0}. If we start with any Kahler metric with potential $\phi$, then use the $a_{\nu}^{(k)}$ as above to  define an algebraic metric with
potential $\phi^{(k)}$ then, after suitable normalisation the $\phi^{(k)}$
converge to $\phi$ as $k\rightarrow \infty$. In particular the algebraic
metrics are dense in the space of all metrics. 

\

\subsubsection{ The Veronese embedding and the Central Limit theorem}

Suppose, in the general situation, that the sections of $L$ generate the
sections of $L^{k}$ so that we have a surjective linear map
$$    s^{k}(H^{0}(L))\rightarrow H^{0}(L^{k}). $$
A metric on $H^{0}(L)$ defines a metric on the symmetric power $s^{k}(H^{0}(L))$
in a standard way. Then we can define a metric on $H^{0}(L^{k})$ by identifying
it with the orthogonal complement of the kernel of the map above. Then we can use this
to define an algebraic Kahler metric on $X$ by the embedding $\iota_{k}: X\rightarrow \bP(H^{0}(L^{k})^{*})$.
Now, up to a  scale factor, these Kahler metrics are independent of $k$. One way of seeing this is that the embedding $\iota_{k}$ is the composite
of $\iota_{1}$£ and the {\it Veronese embedding}
$$   j:\bP(\bC^{N})\rightarrow \bP(s^{k} \bC^{N}), $$
and, up to scale, $j$ is an isometry of the two Fubini-Study metrics.(This
is forced by $U(N)$-invariance.) So the same Kahler metric has a whole series
of algebraic representations.

Let us see how this works in the toric case. We start with data $a_{\nu}$
on $\oP\cap \bZ^{n}$. Then we can write
$$  k \phi= \log \left( \sum a_{\nu} e^{\nu.\ut}\right)^{k}=2 \log \sum B_{\mu}
e^{\mu.\ut}, $$
where the coefficients $B_{\mu}$ are
$$  B_{\mu}= \sum_{\nu_{1}+\dots \nu_{k}=\mu} a_{\nu_{1}}a_{\nu_{2}}\dots
a_{\nu_{k}}.$$
So if we regard $(a_{\mu})$ as a measure $A$ supported on the lattice points
in
 $\oP$ then the $(B_{\mu})$ represent the $k$-fold convolution $A*\dots *A$,
supported on the lattice points in $k\oP$. Now rescale back to the fixed polytope $P$, so we write
$  b_{\nu}^{(k)}= B_{k\nu}$, for $\nu\in \oP \cap k^{-1} \bZ^{n}$. These
define an admissible Kahler potential with Legendre transform $ku$, where
$u$ is the Legendre transform of $\phi$. Then on compact subsets of $P$ we
claim that
\begin{equation}  k^{-1} \log  b_{\nu}^{(k)}= u + O(k^{-1} \log k). \end{equation}

This is essentially the Central Limit theorem, for the convolutions of the
discrete measure $A$. By applying a translation we can reduce to calculating
at the point $\nu=0\in P$. Changing the coefficients $a_{\nu}$ to $a_{\nu}e^{z.\nu}$,
 for any fixed $z\in \bR^{n}$, 
does not change either side of (9), when $\nu=0$, so we can reduce to the case
when $\sum a_{\nu} \nu =0$. That is to say, that $\phi$ attains its minimum
at the point $\ut=0$. Now we consider the function
$$  f(\utheta) = \sum a_{\nu} e^{i \nu. \utheta}. $$
This is a finite trigonometric polynomial which can be regarded as a function
on our compact torus $T$. Then
$$   b_{0}^{(k)}= \int_{T} f^{k} d\utheta , $$
and our assertion follows from the stationary phase approximation, since
the maximum value of $\vert f\vert$ is $  \sum a_{\nu}= u(0)$.

Of course $f$ is just the analytic continuation of $e^{\phi}$, for our Kahler
potential $\phi$. This makes one wonder if there may be other contexts when
it is useful to consider such analytic continuations.

\

\

{\bf Example} For each $k$, the round metric on $S^{2}$ is described as an
algebraic metric with the coefficients $a_{\nu}=  \left( \begin{array}{c} k\\ \nu\end{array}\right)$.

\

\

Notice that the asymptotics approximations we have discussed hold uniformly
over compact subsets of the open polytope $P$. The discussion near the boundary
of $P$ is more delicate, because one gets different asymptotic models. A prototype is the different approximations---normal or Poisson--for the binomial
distribution in different regimes.

   \subsection{Extremal metrics on toric varieties}

   The author has written at length on this topic in other papers, so we
   shall be rather brief here. Expressed in terms of a symplectic potential
   $u$ the condition for an extremal metric is that the scalar curvature
       $$  S(u)=- u^{ij}_{ij}, $$
       is an affine-linear function on $P$. More generally, it is natural
       in this context to consider the prescribed scalar curvature equation
       $ S(u)=A$ for some given function $A$ on $P$. This can be expressed
       as a variational problem. Recall that our polytope $P$ comes with
       preferred defining inequalities $\lambda_{r}(\ux)\geq c_{r}$. These
       linear functions $\lambda_{r}$ define a measure $d\sigma$ on the boundary
       of $P$ (just a multiple of standard Lebesgue measure on each codimension-$1$  face). Then, given a function $A$ on $P$ we define a linear functional
 $$  L_{A}(f)= \int_{\partial P} f d\sigma - \int_{P} A f d\ux. $$
 Now define a nonlinear functional by
 $$  {\cal F}_{A} (u)= L_{A} (u) - \int_{P}\log \det \nabla^{2} u\  d\ux. $$
 
Then an admissible symplectic potential $u$ which satisfies the equation $u^{ij}_{ij}= -A$ is an absolute minimiser of the functional ${\cal F}_{A}$.

The functional ${\cal F}_{A}$ is a variant of the {\it Mabuchi functional},
which is defined in the general Kahler context. It is a convex functional
on the space of convex functions on the polytope $P$. The equation $u^{ij}_{ij}=-A$,
together with the Guillemin boundary conditions asserts that the functional
$L_{A}$ is represented by the inverse of the Hessian of $u$ in the sense
that
\begin{equation}   L_{A}(f)= \int_{P} u^{ij} f_{ij}, \end{equation}
for all test functions $f$. We see immediately from this that if a solution $u$ is
to exist then $L_{A}$ must vanish on the affine linear functions $f$. This
is set of $n+1$ linear constraints on the function $A$. If we take $A$ to
be the constant
$$    \frac{{\rm Vol}(\partial P, d\sigma)}{{\rm Vol}(P, d\ux)}, $$
then $L_{A}$ vanishes on the constant functions $f$. The restriction of this
functional
$L_{A}$ to the linear functions $f$ is the {\it Futaki invariant}, in this
special setting. Otherwise said, this is  essentially the difference between the centre
of mass of $(\partial P, d\sigma)$ in $\bR^{n}$ and the centre of mass of
$(P,d\ux)$. If this Futaki invariant does not vanish then we cannot have
a constant scalar curvature metric, but there is a unique affine-linear function
$A$ satisfying the constraint above, and we seek an extremal metric with this
prescribed scalar curvature.

It is not true that any toric variety admits an extremal metric. To see this
observe that if a solution exists then the weak formulation (10) implies that
$L_{A}(f)\geq 0$ for {\it convex} functions $f$ (with strict inequality if
$f$ is, say, smooth and not affine linear). But one can construct examples
of toric surfaces where $L_{A}$ does not satisfy this condition, for the
affine-linear $A$ above. To fit this in with the discussion of Section 1,
imagine following a minimising sequence $u^{(\alpha)}$ for the functional ${\cal F}_{A}$,
in the case when no solution exists (there would be a similar discussion
for the Calabi functional). Then the typical phenomenon (which one can see
explicitly in some simple examples, and probably holds in general) is that
$u^{(\alpha)}$ behaves like
$$   u^{(\alpha)} \sim C_{\alpha} v, $$
where $ C_{\alpha}$ are real, $C_{\alpha}\rightarrow \infty$ and $v$ is a {\it
piecewise-linear} convex function on $\oP$. Differential geometrically this
corresponds to the {\it collapsing} of some directions in the torus fibration
over the parts of $P$ where the derivative of $v$ is discontinuous. Algebro-geometrically,
the data $v$ describes a toric degeneration of $X$ into a singular toric
variety $X_{0}$ (at least, this is the case if $v$ is defined by \lq\lq rational
data''). In other words we have a picture much like that sketched in 1.1, except
that rather than \lq\lq jumping'' to a different complex structure on the
same underlying smooth manifold we have to allow singularities.  (In fact
a similar thing happens in the Yang-Mills case in higher dimensions, where
the limiting structures may be sheaves rather than holomorphic bundles.)

In this way, one has a good understanding of one mechanism by which existence
can  fail. The
more formidable problem is to see if this is the {\it only} way. More precisely, it is natural to make the
\begin{conj}
If $P\subset \bR^{n}$ is a Delzant polytope  and $A$ is a smooth function
on $\oP$ with the property that $L_{A}(f)$ vanishes if $f$ is affine linear
and $L_{A}(f)>0$ if $f$ is  a convex function 
which is not affine linear, then there is an admissible symplectic potential
satisfying the equation $u_{ij}^{ij}=-A$. 
\end{conj}
We refer to \cite{kn:D4}, \cite{kn:D5}, \cite{kn:D6} for more information about this, particularly
in the case when $n=2$.

  \section{Toric Fano manifolds}
  \subsection{The Kahler-Ricci soliton equation}
  
  The condition that a toric manifold $X$ be Fano, with $L=K_{X}^{-1}$, is
  easily stated in terms of the polytope $P$. There is a preferred \lq\lq
  centre'' $\nu_{0}\in P$ such that for each face $\lambda_{r}(\nu_{0})-c_{r}=1$.
  This follows because the wedge product of the vector fields generating
  the action is a meromorphic $n$-form on $X$ with a simple pole along each
  of the
  divisors corresponding to the faces. Then the inverse is a section of
  $K_{X}^{-1}$ and is a multiple of the standard  basis element $s_{\nu_{0}}$.
  This centre is also the centre of mass of $(\partial P, d\sigma)$.
   
  In this Section we discuss a Theorem of Wang and Zhu \cite{kn:WZ}.
  \begin{thm}
  Any toric Fano manifold has a  Kahler-Ricci soliton metric, unique up to
  holomorphic automorphisms
  \end{thm} 
  
  We will begin by giving a proof which is somewhat different to that of
  Wang and Zhu (although it borrows ideas from that paper and from \cite{kn:TZ}), working largely with
  the symplectic description. We can assume that the centre $\nu_{0}$ is
  the origin. Given a symplectic potential $u$ we write 
  $$  h=x^{i} u_{i} - u, $$
  and $$L=\log \det \nabla^{2} u. $$
  These are smooth functions on $P$ but both tend to infinity at the boundary.
  Note that $h$ depends on a choice of origin in $\bR^{n}$. Of course $h$
  is just the composite of the Kahler potential $\phi$ with the derivative
  of $u$, mapping $P$ to $\bR^{n}$.  The assumption
  that the toric manifold $X$ be Fano is equivalent
  to the fact that, for any admissible  $u$, the difference  $L-h $ is a smooth function on $\oP$. The condition
that $u$ describe a Kahler-Ricci soliton is that 
\begin{equation}   L- h = \sum c_{i} x^{i}, \end{equation}
for constants $c_{i}$ (which of course specify the relevant holomorphic vector
field on the Kahler manifold). Just as in our discussion of extremal metrics,
it is natural in this context to consider more generally an equation
$  L- h= A $ for some prescribed smooth function $A$ on $\oP$. Again, much
as for the extremal case, there are elementary constraints that we need to
impose on $A$. For any symplectic potential $u$ we consider the integrals $$\int_{P} x^{i} e^{L-h} d\ux, $$
for $i=1,\dots,n$. Transforming the integral to the dual space, it becomes
$$  \int_{\bR^{n}} \frac{\partial \phi}{\partial t_{i}} e^{-\phi} d\ut= -\int_{\bR^{n}}
\frac{\partial e^{-\phi}}{\partial t_{i}} = 0.  $$
So a necessary condition that the equation $L-h=A$ has a solution is that,
for each $i$, 
\begin{equation} \int_{P} x^{i} e^{A} d\ux =0  . \end{equation}
This fixes the constants $c_{i}$ in (11). To see this, consider the function of $\underline{c}\in
\bR^{n}$:
   $$ F(\underline{c})= \int_{P} e^{\sum c_{i} x^{i}}d\ux $$
   This is convex and proper (since the origin lies in $P$) and so has a
   unique critical point. But the derivative of $F$ with respect to $c_{i}$
   is 
   $$  \int_{P} x^{i} e^{\sum c_{i} x^{i}} \ d\ux. $$
   So the unique critical point of $F$ gives exactly the constants $c_{i}$
   required to satisfy the constraint.
   
   In sum, the theorem of Wang and Zhu follows from
   \begin{thm}
   For any smooth function $A$ on $\oP$ which satisfies the constraint (12)
   there is a solution $u$ to the equation $L-h=A$, which is unique up to
   the addition of a linear function. \end{thm}
   An equivalent statement is
  \
  
  {\it For any smooth function $A$ on $\oP$ there are  constants $\gamma_{i}$
  and an admissible potential $u$ such that $L-h=A+\sum \gamma_{i} x^{i}$.
  The $\gamma_{i}$ are unique and $u$ is unique up to the addition of a linear
  function.}
  \
  
  The equivalence of the statements follows from the same argument as above.

\subsection{Continuity method, convexity and a fundamental inequality}

For any symplectic potential $u$ on our Fano polytope, centred at the origin,
we write $\rho=L-h$. Now we define the following weighted norms, for
functions $f,g$ on $P$:

$$   \langle f, g\rangle_{u}= \int_{P} f g e^{\rho}\  d\ux; $$
$$  \langle \nabla f, \nabla g \rangle_{u}= \int_{P} f_{i} g_{a} u^{ia} e^{\rho}\
d\ux;$$
$$ \langle \nabla^{2} f , \nabla^{2} g\rangle_{u}=\int_{P} f_{ij} g_{ab}
u^{ia} u^{jb} e^{\rho}\ d\ux. $$

The first variation of $\rho$ with respect to an infinitesimal variation
$f$ in $u$ is
$ \delta \rho= \Box f$,
where $\Box$ is the differential operator
\begin{equation} \Box f = u^{ij} f_{ij}- x^{i} f_{i}+ f. \end{equation}
Since $$ \rho_{j}= - u^{ia}_{a} u^{aj} - x^{a} u_{ja}, $$

this can also be written as 
\begin{equation} \Box f = \left(u^{ij} f_{i}\right)_{j} - u{ij} \rho_{j} f_{i} + f, \end{equation}
from which it follows that
\begin{equation}   \langle \Box f, g\rangle_{u}= -\langle \nabla f ,\nabla g \rangle_{u}+
\langle f, g \rangle_{u}. \end{equation}
In particular, $\Box$ is self-adjoint with respect to the weighted norm.

Now define a functional by
\begin{equation} {\cal F}(u)= \int_{P} e^{\rho}\ d\ux. \end{equation}
Then the first variation is
\begin{equation} \delta {\cal F}= \int_{P} \Box f e^{\rho} d\ux = \langle \Box f, 1\rangle_{u}. \end{equation}
By the self-adjoint property we can also write this as
\begin{equation} \delta {\cal F}= \langle f, \Box 1\rangle_{u} = \langle f, 1\rangle_{u}=\int_{P}
f e^{\rho} d\ux. \end{equation}
This leads to two different expressions for the second variation of ${\cal
F}$. If we put $u_{t}= u + t f, \rho_{t}= \rho(u_{t})$ and write $\Box_{t}$ for the operator defined
by $u_{t}$ then
$$  \frac{d}{dt} \Box_{t} f = - u^{ia} u^{jb} f_{ij}f_{ab}. $$
So, 
$$  \frac{d^{2}}{dt^{2}} {\cal F}(u_{t})= \frac{d}{dt} \int_{P} \Box_{t}
f e^{\rho_{t}} d\ux= \int_{P} \left( \Box_{t}f \Box_{t} f - u^{ia}u^{jb} f_{ij} f_{ab}
\right)e^{\rho_{t}} d\ux, $$
which is equal to $$ \langle \Box_{t} f, \Box_{t} f \rangle_{u_{t}}- \langle
\nabla^{2} f , \nabla^{2} f \rangle_{u_{t}}. $$
On the other hand
$$  \frac{d^{2}}{dt^{2}}{\cal F}(u_{t})= \frac{d}{dt} \int_{P} f e^{\rho_{t}}
d\ux= \int_{P} f \Box_{t} f e^{\rho_{t}} d\ux. $$
So, evaluating at $t=0$ and dropping $t$ from the notation, we have the identity
\begin{equation} \langle \Box f, \Box f \rangle_{u}- \langle \nabla^{2} f, \nabla^{2}
f \rangle_{u} = \langle f, \Box f \rangle_{u}. \end{equation}
Applying (15), with $g=\Box f$, this gives,
\begin{equation}  \langle \nabla f, \nabla \Box f\rangle_{u}= -\langle \nabla^{2} f, \nabla^{2}
f \rangle_{u}. \end{equation}
It is obvious from the definition that $\Box$ vanishes on the linear functions
and $\Box 1=1$. If $f$ is any eigenfunction of $\Box$, with eigenvalue $\lambda$, which is orthogonal
to the linear functions and then constants, then $\nabla^{2} f$ is non-zero
and the identity gives
$$  \lambda \langle \nabla f, \nabla f \rangle_{u}= -\langle \nabla^{2} f,
\nabla^{2} f \rangle_{u}, $$
so $\lambda<0$.  (This is a variant of the standard lower bound on the eigenvalues
of the Laplacian on a manifold with positive Ricci curvature, the identity
can of course be verified more directly, but the argument above avoids some
laborious  manipulation.) In sum, we have derived an inequality
\begin{equation}   \langle 1, f\rangle_{u}\langle 1,1\rangle_{u} -\langle f, \Box f\rangle_{u} \geq 0, \end{equation}
with equality if and only if $f$ is a linear function.

\
Now to apply this to our problem. First, we can use the continuity method
for the equation $L-h=A+\sum \gamma_{i} x^{i}$, with respect to variations in $A$. The linearised
equation is $\Box_{u}f=\delta A+\sum \delta \gamma_{i} x^{i}$. Since the
cokernel of $\Box_{u}$ is identified with the linear functions this linearised
equation has a solution and we can apply the implicit function theorem in
the usual way.

Second, we obtain the uniqueness of solutions. Consider the functional
$-\log {\cal F}$. Along a line $u_{t}= u+tf$ we have
$$   \frac{d^{2}}{dt^{2}}\left( -\log{\cal F}\right)= \frac{1}{{\cal F}^{2}}( {\cal F}{\cal
F}'- {\cal F}'')$$
where ${\cal F}', {\cal F}''$ denote the derivatives of ${\cal F}$. Evaluating
at $t=0$ we have
$$ {\cal F}=\langle 1,1\rangle_{u}, {\cal F}'=\langle 1,f\rangle_{u}, {\cal F}''=
\langle f, \Box_{u} f\rangle_{u}, $$
so our inequality (21) asserts that the second derivative of $-\log {\cal F}$ is positive, and strictly positive unless $f$ is affine-linear. Thus $-\log{\cal
F}$ is a convex function. Now if $\rho=A+\sum \gamma_{i}x^{i}$ the $\gamma_{i}$
are determined by $A$, using the same argument as in the previous subsection.
So we may as well suppose that $\gamma_{i}=0$. Then ${\cal F}(u)= C$ where
$C$ is the integral of $e^{A}$. The equation $\rho=A$ is the Euler-Lagrange
equation for critical points of the linear function
$$ u\mapsto \int_{P} u e^{A} \ d\ux $$
subject to the constraint $-\log {\cal F}= -\log C$. The convexity gives
uniqueness, modulo linear functions.

   \subsection{A priori estimate}
   
   To prove Theorem 2 we need to establish appropriate {\it a priori}
   bounds on a solution to our equation. We proceed in four steps.
   
   \newpage
   
   {\it Step 1: Preliminaries}
   
   \
   
   We want to appeal to some of the standard body of theory for compact Kahler manifolds,
   that is, where we consider a fixed reference metric $\omega_{0}$ on a compact
   manifold and another metric $\omega=\omega_{0}= i \dbd \psi$.
       Our problem differs a little from that usually considered in
       the literature. To fit into a general  setting we could consider a
       fixed smooth function $G$ of $p$-variables, a compact Kahler manifold
       $X$ with $p$ fixed  holomorphic vector fields $v_{\alpha}$  and a function $\psi$
       which satisfies an equation
       $$    (\omega_{0} + i \dbd \psi)^{n} = \exp(\psi + G(\nabla_{1} \psi,
       \dots , \nabla_{p} \psi)) $$
       where $\nabla_{\alpha} \psi$ denotes the derivative of $\psi$ along
       the vector field $v_{\alpha}$. Then the modification by  Tian and Zhu (\cite{kn:TZ}, Section 5, especially Prop. 5.1)
       of the standard argument of Yau, shows that in this situation
       an $L^{\infty}$ bound  on $\psi$ leads to bounds on all higher derivatives.
       (Apart from this the proof we give is self-contained.)
       
        In our toric setting, we choose
   some fixed admissible Kahler potential $\phi_{0}$ on $\bR^{n}$ with Legendre
   transform $u_{0}$. Then we consider some general Kahler potential $\phi$,
   with Legendre transform $u$ and set $\psi=\phi-\phi_{0}$. So an $L^{\infty}$
   bound on $\psi$ on the compact toric manifold is identical to an $L^{\infty}$
   bound on $\phi-\phi_{0}$ on $\bR^{n}$. Now a general property of the Legendre
   transform is that it is an {\it isometry} with respect to the $L^{\infty}$ distance:
   that is to say
   $$  \sup_{\ut\in \bR^{n}} \vert \phi(\ut)-\phi_{0}(\ut)\vert= \sup_{x\in
   P} \vert u(x)- u_{0}(x) \vert . $$
    
   This is an elementary exercise.
  
       \
       
       In our situation, $u_{0}$ is a fixed continuous function on $\oP$
       so an $L^{\infty}$ bound on the function $\psi$ on the compact Kahler
       manifold is equivalent to an $L^{\infty}$ bound on the \lq\lq unknown''
       symplectic potential $u$.

       In sum,  we see that to  prove our
       proposition it suffices to establish an {\it a priori} $L^{\infty}$ bound on 
       symplectic potentials $u$ satisfying a differential inequality
          \begin{equation}\vert L- h \vert \leq C, \end{equation}
          for fixed $C$. 
          Of course for this to make sense we have to normalise the non-uniqueness
          under the addition of linear functions, but we can do this very         simply     by
          restricting to functions $u$ whose derivative vanishes at the origin.
          i.e are minimised at the origin. 
          We write $m=-u(0)$,  so  our problem comes down
          to obtaining upper and lower bounds on $m$ and an upper bound on
            $\rm {Max}_{\oP}
          u - u(0)$. By the Sobolev inequality, the latter will follow if
          we establish an {\it a priori} bound on $\Vert \nabla u \Vert_{L^{p}}$,
          for any $p>n$.

          \
          
          {\it Step 2}
          \
          
          Here we  get
          a lower bound on $m=-u(0)$. By definition $h(0)=m$. Let the polytope $P$ be contained in
          the $R_{1}$- ball about $0$ in $\bR^{n}$ and  let $\Omega\subset P$ be the set where $\vert \nabla u \vert\leq
          1$. Since the derivative of the Legendre transform $\phi$ is bounded
          in size by $R_{1}$ we have  $\vert h(x)-m\vert\leq R_{1}$  for
          all $x\in \Omega$
          and the basic assumption (22) gives
          $   L\leq m+R_{1}+ C$ so
          $$  \det(\nabla^{2} u)\leq \exp(m+R_{1}+C). $$
          But the integral of $\det(\nabla^{2} u)$ over $\Omega$ gives the volume
          $\omega_{n}$ of the unit ball in $\bR^{n}$ so
          $$  \exp(m+R_{1}+C) \Vol(\Omega)\geq \omega_{n}. $$
          Since the volume of $\Omega$ cannot exceed the volume of $P$ this
          gives a lower bound on $m$.
          
          \
          
          {\it Step 3}
          
          \
          
          Here we bring in a crucial identity. 
          \begin{lem}
          For any admissible symplectic potential $u$ and each index $i$ we have
          $$  \int_{P} e^{\rho} u_{i} x^{i} \ d\ux =  \int_{P} e^{\rho}
          \ d\ux , $$
          where $\rho=L-h$.
          \end{lem}
          To see this,  transform the integral to $\bR^{n}$, giving
          $$\int_{P}e^{\rho} u_{i} x^{i}\ d\ux = \int_{\bR^{n}} e^{-\phi}
          \frac{\partial \phi}{\partial t_{i}} t_{i} \ d\ut. $$
          Now integrate by parts to write this as
          $$  \int_{\bR^{n}} e^{-\phi} d\ut, $$
          which transforms back  to
          $$  \int_{P} e^{\rho}. $$
          
          \
          
          \
          
          In our situation we are given $\vert \rho \vert\leq C$, so
          summing over the index $i$, we have
          $$  \int_{P} \sum u_{i} x^{i}  d\ux \leq n  e^{2C} \Vol(P). $$
          In terms of \lq\lq generalised polar co-ordinates'' $(r,\utheta)$
       we can write this as
       $$  \int_{P}  r^{n} \frac{\partial u}{\partial r}  dr d\utheta \leq
       n e^{2C} \Vol(P). $$
          By hypothesis, the derivative of $u$ vanishes at the origin and
          the radial derivative $\frac{\partial u}{\partial r}$ is increasing
          along each ray. It is clear then that if we fix some
           $R_{0}$ such
          that the $3R_{0}$- ball $B_(3 R_{0})$ about the origin is contained in $P$ then
          we get from Lemma 1 an {\it a priori} bound on the integral of $u-u(0)$ over the
          boundary of $B(2 R_{0})$. In turn this gives, by an easy elementary
          argument, a bound on the derivative over an interior ball, say
          \begin{equation} \vert \nabla u \vert \leq C' \ \ {\rm on} \ \ B(R_{0}) , \end{equation}
          where $C'$ can be computed explicitly in terms of $n,C, R_{0}, \Vol(P)$.
          
          \
          
          \
          
          {\it Step 4}

          \
          
          We can now complete our task,  using some elementary convexity
          arguments. First, throughout $P$ we have $h\geq m$ so $\det\nabla^{2}
          u \geq e^{-C} e^{m}$. Thus the  volume of the image of the ball
          $B(R_{0})$
          under the map $Du$ is at least $e^{-C} e^{m} \Vol(B(R_{0}))$. But
         by (23) this image is contained in the ball of radius $C'$, so
         $$  e^{-C} e^{m} \leq \left( \frac{C'}{R_{0}}\right)^{n}. $$
         This gives our upper bound on $m$. 
         
         Next we claim that at any point $x_{0}\in P$ where $h(x_{0})=m+1$ we have
         $\vert \nabla u \vert \leq \kappa$, where $\kappa=C'+R_{0}^{-1}$. To see this, consider the affine-linear function $\pi$ defining
         the supporting hyperplane of $u$ at $x_{0}$. By definition, this
         has the form $\pi(x)= -h(x_{0}) + D(x)= -(m+1) + D(x) $ where $D$
         is a linear function, equal to the derivative of $u$ at $x_{0}$.
         Thus there is a point $y$ on the boundary of $B(R_{0})$ where
         $\pi(y)= -(m+1) + R_{0} \vert D \vert$. By convexity we have
         $u(y)\geq \pi(y)$, so $u(y) \geq -(m+1) + R_{0} \vert D\vert$. On
         the other hand (23) implies that $u(y)\leq - m + C' R_{0}$, so we have
         $$   -(m+1) + R_{0} \vert D \vert \leq - m + C' R_{0}, $$
         and $\vert D\vert \leq C' + R_{0}^{-1}$.
         
         Now consider the situation under the Legendre transform. The statement
         above becomes that for any point $\ut_{0}$ where $\phi(\ut_{0})=\phi(0)+1=m+1$
         we have $\vert \ut_{0} \vert \leq \kappa$. It follows easily from the convexity
         of $\phi$ that for {\it any} $\ut$ we have
         \begin{equation}\phi(\ut) \geq \kappa^{-1} \vert \ut \vert + m -
         \kappa^{-1}. \end{equation}
         Recall from the conclusion of Step 1 that we seek a bound on $\Vert
         \nabla u \Vert_{L^{p}}$. By the definition of the Legendre transform
         we have
         $$  \int_{P} \vert \nabla u \vert^{p} \ d\ux = \int_{\bR^{n}} \vert
         \ut \vert^{p} \det (\nabla^{2} \phi)\ d\ut. $$
         Using our hypothesis (22) we have then
         $$ \Vert \nabla u \Vert_{L^{p}}^{p} \leq e^{C} \int_{\bR^{n}}  \vert \ut
         \vert^{p} e^{-\phi} \ \ d\ut. $$
         Now write
         $$ \int_{\bR^{n}} \vert \ut \vert^{p} e^{-\phi} \ d\ut \leq e^{m-\kappa^{-1}}
         \int_{\bR^{n}} \vert \ut\vert^{p} e^{- \kappa^{-1} \vert \ut \vert}\
          d\ut, $$
          to finish the argument.
         
         \
         
         \
         
       {\bf Remark}
     
     The identity in Lemma 1 is related to the positivity condition for the
     linear functional $L_{A}$ discussed in Section
    2. Suppose we have a Fano polytope with vanishing Futaki invariant, so
    we seek a Kahler-Einstein metric which {\it a fortiori} has constant
    scalar curvature. Thus, in the setting of Section 2, the function $A$
    is a constant. By an argument of Zhou and Zhu \cite{kn:ZZ} the linear
    functional $L_{A}$ satisfies the positivity condition in this case. On
    the other hand, in the setting of Section 2, anytime $L_{A}$ satisfies
    the positivity condition we can deduce an {\it a priori} interior  bound on the
    derivative of a normalised solution, see \cite{kn:D4}, just as in (23)
    above. Thus the identity of Lemma 1 can be seen as an extension of this
    discussion to the case of general Ricci soliton metrics. 
    
    \
    
    \

  \subsection{The method of Wang and Zhu}

  We will now discuss briefly the original approach of Wang and Zhu (again
  with some modifications). For simplicity
  we will just consider the case when the Futaki invariant vanishes, so we
  seek a Kahler-Einstein metric. Recall from the above that the vanishing
  Futaki invariant is equivalent to fact that the centre of mass of the polytope
  $P$ is the  preferred centre, which we are taking as $0\in \bR^{n}$. 
  
  Wang and Zhu
   use the continuity method with respect to the family of equations
       
       \begin{equation}  \det(\nabla^{2} \phi) = \exp(-(s\phi+(1-s)f)), \end{equation} where $f$ is a fixed
       admissible Kahler potential and $0\leq s<1$. 
       We  discuss first the case when $s=0$. Then the equation in question
       is just the toric case of the \lq\lq prescribed volume form'' equation, solved, for general Kahler manifolds, by Yau. But let us see how to give a simple proof in this special situation.
      The equation (25) with $s=0$ is degenerate, in that we can obviously
       change $\phi$ by the addition of a constant, so we may normalise $u$
       to be zero at some point. Thus, as before, all we need to do is bound the $L^{p}$
       norm of $\nabla u$. But for this we simply write
         
$$  \int_{P} \vert \nabla u \vert^{p}\ d\ux = \int_{\bR^{n}} \vert \ut\vert^{p}
\det \nabla^{2} \phi\  d\ut =\int_{\bR^{n}}\vert \ut \vert^{p} e^{-f}\ d\ut<\infty.
$$ This concludes the proof of the $L^{\infty}$ estimate for the case $s=0$. (Here we have
not used the fact that $f$ is convex, so by deforming $f$ one can prove the
toric case of  Yau's
Theorem: the existence of a solution for any $f$.) 

Now we go on to the main case, when $s>0$. It suffices to obtain estimates
for $s\geq s_{0}$ for some fixed $s_{0}>0$.

        Set $w=s\phi+(1-s)f$. Then $w$ is another admissible function
       and 
       $$  \det(\nabla^{2} w)\geq s_{0}^{n} \det (\nabla^{2} \phi) = s_{0}^{n} e^{-w}.
       $$
    
       Let the minimal value of $w$ be $m$, attained at a point $\zeta\in
       \bR^{n}$.
       The first main step in the proof is
       \begin{prop}
       We have
       $$w(\ut) \geq \epsilon \vert \ut-\zeta \vert - C $$
       for known $\epsilon, C$.
       \end{prop}

       The foundation of the approach of Wang and Zhu is the following fact.
       
      \begin{prop}
      
      Suppose that $v$ is a  convex function on $\bR^{n}$, attaining minimal value
      $0$, and suppose $\det(\nabla^{2} v)\geq
      \lambda$ when $v\leq 1$. Then if $K$ is the set where $v\leq 1$ we
      have $\Vol(K)\leq C \lambda^{-1/2}$ for some constant $C$ depending
      only on the dimension $n$.
      \end{prop}
      
      Wang and Zhu prove this using a  comparison argument. It can also be shown using
      the elementary geometry of the derivative of $v$ (see \cite{kn:Gut}
      Prop. 3.2.3), but both
      approaches depend on the fact that after a unimodular affine transformation
      we can suppose that there are concentric balls
      $$   B(R_{1}) \subset K \subset B(R_{2}), $$
      with the ratio $R_{2}/R_{1}$ of the radii bounded by a fixed constant
      depending on the dimension. Notice that a reverse inequality holds.
      If in the same situation $\det(\nabla^{2} v) \leq \Lambda$ then
      $\Vol(K)\geq C\Lambda^{-1/2}$ (\cite{kn:Gut}, Cor. 3.2.4).
      
      \

      With this background in place we can proceed to explain the proof of Wang
      and Zhu.
      Let $m$ be the minimal value of the function $w$ and set $v=w-m$.
      Then $\det(\nabla^{2} v)\geq \lambda= t_{0}^{n} e^{m+1}$ on the set $K$ where
      $v\leq 1$. So we deduce that
      \begin{equation}  \Vol(K)\leq C \lambda^{-1/2}= C' e^{m/2}, \end{equation}
      say.  For each positive $\mu$ let $K_{\mu}$ be the set $\{ v\leq \mu \}$ and
      $V(\mu)=\Vol(K_{\mu})$. Then
      convexity implies that $K_{\mu}$ is contained in the dilate of $K$ by
      factor $\mu$ about the minimum point of $v$. Thus
      $$  V(\mu)= \Vol(K_{\mu})\leq h^{n} \Vol(K) \leq \mu^{n} C' e^{m/2}. $$
      By the co-area formula
           $$  \int_{\bR^{n}} e^{-w} \ d\ut= \int_{0}^{\infty} e^{-\mu} V(\mu)
           \ dh. $$
      Now the volume form $\det(\nabla^{2} \phi)$ is at most $ e^{-m} e^{-v}$
      and its integral is the volume of our manifold
      $X$.
      So
$$  \Vol(X) \leq e^{-m}  \int_{0}^{\infty} C' e^{m/2} e^{-\mu} \mu^{n}\ d\mu=
C'' e^{-m/2}, $$
say. We see  that
  \begin{equation} m\leq m_{0} = 2 \log(I_{0}/C''),  \end{equation}
  and then deduce from (26) that
  \begin{equation}\Vol(K) \leq C' e^{m_{0}/2}. \end{equation}
  
  Now we use the fact that $\vert \nabla w\vert \leq b$ say. This means
  that the distance from the boundary of $K$ to the minimum point $\zeta$.
  is at least
  $b^{-1}$, so $K$ contains a ball of this fixed  radius about $\zeta$. If
   $K$ contains a point $\zeta'$ with $\vert \zeta-\zeta'\vert =R$
  for large $R$, then the volume of $K$ would be large, contradicting the
  bound (28). So
  we conclude that $K$  is contained in the ball $\{\zeta':\vert \zeta'-\zeta\vert
  \leq R_{0}\}$ for some fixed $R_{0}$. But then convexity implies that
  $$  \vert \xi-\zeta \vert \leq R_{0}^{-1} v(\xi). $$
  This completes the proof of   Proposition 1. 
  
  \
  
  The second main step is to show that $\vert \zeta\vert $ is not large.   This is where the hypothesis that the  the Futaki invariant vanishes is
  used.  Consider the
  derivative $Df$ of the fixed admissible function $f$. This is a vector-valued function on $\bR^{n}$, which gives
  a proper map to the open polytope $P$. 
  The crucial thing is an identity
  \begin{equation} \int_{\bR^{n}} Df e^{-w} \ d\ut= 0. \end{equation}
  To see this, consider one component $\frac{\partial f}{\partial t_{a}}=
  f^{a}$ of $Df$, and observe first that
  $$   \int_{\bR^{n}} \left((1-s) \frac{\partial f}{\partial t_{a}}+s \frac{\partial \phi}{\partial t_{a}}\right) 
  e^{-w} \ d\ut= \int_{\bR^{n}} \frac{\partial w}{\partial t_{a}} e^{-w}
  \ d\ut= 0.
  $$
  So it is the same to show that
  $$ \int_{\bR^{n}} \frac{\partial \phi}{\partial t_{a}} e^{-w} \ d\ut =0.$$
  But this integral is
  $$  \int_\bR^{n} \frac{\partial \phi}{\partial t_{a}}  \det(\phi_{ab})\
  d\ut $$
  which is the same as
  $$  \int_{P} x^{a}\ d\ux $$
  and this vanishes by our hypothesis.
  
  \
  
  Consider a codimension-$1$ face of $P$ defined by an equation $\lambda_{r}(x)=c_{r}$.
  Let $g_{r}$ be the function 
  $$  g_{r}(\ut)= \log ( \lambda_{r}(Df(\ut))- c_{r}). $$
  It is easy to check that the derivative of $g_{r}$ is {\it bounded} on $\bR^{n}$. Suppose
  $\vert \zeta\vert$ is large. This means that $Df(\zeta)$ is close to the
  boundary of $P$, so there is some $r$ for which $g_{r}(\zeta)$ is very negative
  $  g_{r}(\zeta)\leq -M$ say, for $M$ large. Then the bound on the derivative
  of $g_{r}$ means that we can find a constant $\sigma$ such that on the
  ball $B$ of radius $\sigma M$ about $\zeta$ we  have $g_{r}\leq -M/2$.
  Thus $\lambda_{r}(Df)\geq c_{r}/2$ say, on $B$, if $M$ is large enough. Equally,
  it follows from Proposition 1 that when $M$ is  large the integral of $e^{-w}$ over
  $\bR^{n}\setminus B$ is small. This shows that 
  $$\int_{\bR^{n}} \lambda_{r} (DF)  e^{-w} >0, $$
  if $M$ is large, which is a contradiction to the identity (29)  above.

  It is now easy to complete the proof. Since $\zeta$ is bounded we have
  $$  w(\ut)\geq \epsilon \vert \ut \vert - c $$
  and the bound on the $L^{p}$ norm of $\nabla u$ follows just as before. Then it is straightforward to get upper and lower bounds on $u$
  at some point, for example
 the point corresponding to $\zeta$.

      \section{Variants of toric differential geometry}
      \subsection{Multiplicity-free manifolds}
      
      The  special features of toric differential geometry can be traced
      back to the fact that the group of Hamiltonian diffeomorphisms which commute
      with the  action is abelian. In general, the action of a compact group $G$
       on  a symplectic manifold $(M,\omega)$ is called \lq\lq multiplicity-free''
       if it has this property. This is equivalent to saying that the all
       the $G$-invariant functions Poisson-commute. The theory has been developed
       by a number of authors. The analogous notion in algebraic geometry
       is that of a {\it spherical variety}. The theory of extremal metrics
       and the Mabuchi functional in this setting
       has been studied by Alexeev and Katzarkov \cite{kn:AK} and by Raza \cite{kn:Raz}
       and 
       Podesta and Spiro \cite{kn:PS1} have extended the theorem of Wang and Zhu for Fano manifolds in this direction. There is also related
work of Bielwaski \cite{kn:Biel}. We will now outline some of these
ideas.

       There is  a general classification of multiplicity-free manifolds
       (\cite{kn:Wood}, \cite{kn:L}),  but rather than attempting to discuss the
       most general situation we focus on a simple class of examples.
       Pick a maximal torus  $T$ in the compact connected Lie group $G$ and let
       $V$ be the dual of the Lie algebra of $T$. There is a weight lattice
       $\Lambda \subset V$. Pick a positive Weyl chamber in $V$ and consider
       an integral Delzant polytope $P$ whose closure is contained in the
       interior of this chamber. We construct a manifold from this data and
       as usual we can take either a symplectic or complex point of view.
       
       {\it Complex}
       
       The choice of a Weyl chamber defines a Borel subgroup $B$ of the complexified
       group $G^{c}$, containing
       the complexified torus $T_{c}$. For example if $G=U(m)$ the Borel
       subgroup is the group of complex matrices with zeros below the diagonal.
       Then we have a generalised flag manifold $Y=G^{c}/B$, which is a compact
       complex manifold. There is a homomorphism from $B$ to $T^{c}$ which
       is a left inverse to the inclusion. Now form the toric manifold $X$
       associated to the polytope $P$. Then $T^{c}$ acts holomorphically on $X$ and so
       $B$ does also via the homorphism above. So we get a complex manifold
       \begin{equation}  Z=  G^{c}\times_{B} X, \end{equation}
       with a holomorphic fibration $\pi:Z\rightarrow Y$,  having fibre $X$.
       The group $G^{c}$ acts on $Z$ and $\pi$ is a $K^{c}$-equivariant
       map.  Further,
       we have a $T^{c}$-equivariant line bundle $ L\rightarrow X$ so the
       same construction yields a $G^{c}$-equivariant line bundle $\cL\rightarrow Z$ which restricts
       to $L$ on each fibre. We can identify $H^{0}(Z,\cL)$ with the sections
       of the vector bundle $\pi_{*}(\cL)$ over $F$. Recall that there is
        a standard basis for $H^{0}(X,L)$ labelled by the lattice
       points $\nu$ in $\oP$. This yields an isomorphism between $\pi_{*}(\cL)$
       and the direct sum of line bundles $\xi_{\nu}\rightarrow F $ associated
       to these weights. The Borel-Weil theorem asserts that the holomorphic
       sections of $\xi_{\nu}$ define the irreducible representation $W_{\nu}$
       of $G^{c}$ with highest weight $\nu$. So we see that, as a representation
       of $G^{c}$,
       $$   H^{0}(Z,\cL)= \bigoplus_{\nu\in \oP} W_{\nu}. $$
       In particular the representation is \lq\lq multiplicity-free'', in
       the sense that all irreducibles appear with multiplicity at most one.
       This is the same as saying that the algebra of $G^{c}$-equivariant
       endomorphisms of $H^{0}(Z,\cL)$ is commutative. The terminology \lq\lq
       multiplicity free'' in the symplectic setting is derived by analogy
       with this.
       
       Notice that replacing $P$ by a multiple $kP$ yields the same complex
       manifold $Z$ but replaces $\cL$ by $\cL^{k}$. Translating $P$ by $\nu$
       does not change $Z$ but 
       changes the line bundle $\cL$ to $\cL\otimes \pi^{*}(\xi_{\nu})$.
       In
        none of the above do we use the fact
       that $\oP$ lies in the interior of the positive Weyl chamber. This
       is exactly the condition which implies that $\cL$ is an {\it ample}
       line bundle over $Z$. 
       
       \
       
       {\bf Example}
       Take $G=SU(2)$, so $V$ can be identified with $\bR$ and the positive
       Weyl chamber with the positive reals. Let $P$ be the interval
       $(p_{1}, p_{2})$. Then $X=Y=\bC\bP^{1}$ and $Z$ is the blow-up of the
       complex projective plane atone point. As $p_{1}, p_{2}$ vary we get all positive
       line bundles $\cL$
       over $Z$.

       \

       For the symplectic description we start by writing $Y=G/T$, and think
       of $G$ as a principal $T$-bundle over $Y$. As a manifold $Z$ is the
       associated bundle $G\times_{T} X$.  Now $T$ has a Hamiltonian
       action on $X$. In general suppose a Lie group $K$ has a Hamiltonian action
       on a symplectic manifold $(M,\Omega)$ and we have a  principal
       $K$-bundle $E\rightarrow U$. Then there is a canonical closed $2$-form
       $\tilde{\Omega}$ on the associated bundle $E\times_{K} M$ which restricts to $\Omega$
       (in the obvious sense) on each fibre. Indeed this is true in the \lq\lq
       universal'' case when we take the group of all Hamiltonian diffeomorphisms
       of a symplectic manifold. This theory is explained in detail in \cite{kn:McDSal},
       Sect. 6.1).
       It is easy to say explicitly how this works in the case at hand. Choose a basis of $V=\Lie(T)^{*}$.
The basis elements can be regarded as left-invariant $1$-forms $\alpha_{i}$
on $G$ and  also as the components of a connection form on the $T$-bundle
$G\rightarrow Y$. The moment map $\mu:X\rightarrow V$ has components, relative
to this basis, which we denote by $x^{i}$, in line with our previous notation.
Since the moment map is equivariant we can also regard $\mu$ as a map from
$Z$ to $V$ and the components $x^{i}$ as functions
on  $Z$. Restrict to the open set $Z_{0}\subset Z$ corresponding
to the open set $X_{0}\subset X$ where $T$ acts freely. This can be identified
with the product $ P\times G$, so we can also regard $\alpha_{i}$ as $1$-forms
on $X_{0}$. Then we set
$$  \tilde{\Omega} = d( \sum x^{i} \alpha_{i}) $$
on $Z_{0}$. On each fibre the $1$-forms $\alpha_{i}$ can be identified with
the $d\theta_{i}$ and we recover the form $\sum dx^{i}d\theta_{i}$. The point
is that, although the $1$-forms $\alpha_{i}$ do not extend over $Z$, the closed
$2$-form $\tilde{\Omega}$ does. This is fairly clear from the corresponding
discussion
on the fibres. The condition that $\oP$ lies inside an open Weyl chamber
is exactly the condition that the form $\tilde{\Omega}$ is symplectic. The
$G$-invariant functions on $Z$ are just the composite of $\mu$ with functions
on $\oP$ and these all Poisson-commute.

 An important object in this theory 
 is the \lq\lq Duistermaat-Heckmann''function $W$ on $V=\Lie(T)^{*}$. It is
 a polynomial function which, on the  open Weyl chamber, gives the symplectic
 volume of the corresponding coadjoint orbit. Algebraically it is the product
 of the positive roots, where the roots are viewed as linear functions on
 $V$.  The push-forward $\mu_{*}(\tilde{\Omega}^{N})$of the  symplectic measure on $Z$ is the restriction to $\oP$ of $(2\pi)^{n} W$ times the Lebesgue measure on $V$. Thus if we identify functions on $\oP$ with $G$-invariant
functions on $Z$ the operation of integration over $Z$ corresponds to  the
weighted integral
    \begin{equation} \int_{\oP}  f W d\ux. \end{equation}

Raza extended the symplectic point of view on  toric differential geometry,
as outlined (2.1.2) above,   to this
setting \cite{kn:Raz}. The orthogonal complement with respect to $\tilde{\Omega}$ defines
a field of horizontal subspaces in $Z$, transverse to the fibres. Any $G$-invariant
almost-complex structure on $Z$, compatible with $\tilde{\Omega}$, must respect
this decomposition and agree with the standard complex structure, induced
from $Y$, in the horizontal subspace. So such almost-complex structures correspond
to the same $T$-invariant almost-complex structures on $X$ which we studied
before, and the integrable structures are determined by an admissible symplectic
potential $u$ on $\oP$, as before. The whole difference in the theory resides
in the weight function $W$. Raza shows that the scalar curvature of the metric
on $Z$ defined by a symplectic potential $u$ is
$$       \frac{1}{W} \ \frac{\partial^{2} W u^{ij}}{\partial x^{i} \partial
x^{j}} + f_{G}, $$
where $f_{G}$ is function determined by the group $G$. In fact if we let
$\sigma\in \Lie(T)^{*}$ be the sum of the positive roots of $G$ then
$$   f_{G}= W^{-1} ( W_{i} \sigma^{i}): $$
the derivative  of $\log W$ in the direction $\sigma$. This extends Abreu's
formula in the toric case, and also a formula of Calabi, for the case when
$K=SU(2)$ (\cite{kn:Cal}, \cite{kn:HS}). There there seems to be considerable scope for extending the analytical theory developed in the toric case to this more
general setting, similar to the work of  Szekelyhidi in \cite{kn:Sz}.

\

Now we consider the Fano case, where the line bundle ${\cal L}$ is $K_{Z}^{-1}$.
This requires, first, that the fibre $X$ be Fano. Recall that there is a
preferred centre $\nu_{0}$ in $P$ (the centre of mass of the boundary). The second requirement,
to identify ${\cal L}$ with $K_{Z}^{-1}$, is that $\nu_{0}$ is equal to $\sigma$,
the sum of the positive roots. (To see this, observe that the line bundle
over $Y$ associated to the weight $\sigma$ is the $K_{Y}^{-1}$.)   In Section
3 we took this centre to be the origin, but here that would conflict with the Weyl chamber
structure. So, given a polytope $P$ satisfying these two conditions above,
and an
admissible symplectic potential $u$, we define
$$  h= (x^{i}- \sigma^{i}) u_{i} - u. $$
Then  $L-h$ is smooth on $\oP$. The Ricci soliton  condition is
$$ L-h= G+ \sum c_{i} x^{i}, $$
for suitable constants $c_{i}$. This falls into the class of equations we considered
in 3.2, and the existence theorem of Podesta and Spiro is another illustration
of our result there.

What we have discussed is the simplest class of multiplicity-free manifolds.
One gets other examples in at least two ways.
\begin{itemize}
\item One can allow the boundary of $\oP$ to touch the boundary of the Weyl
chamber. 
\item One can consider polytopes contained in proper affine subspaces of
$\Lie(T)^{*}$. 
\end{itemize}
There seems to be considerable scope for developing this theory, both in
the Fano case and for extremal metrics. In the latter case one could hope
to extend the results proved for toric varieties, along the lines of the
work of Szekelyhidi \cite{kn:Sz} in the case when $G=SU(2)$.

      \subsection{Manifolds with a dense orbit}
      
      Now we consider another generalisation of toric geometry. Let $G$ be
      a compact Lie group and $G^{c}$ its complexification. Suppose $G^{c}$
      acts holomorphically on a compact complex manifold $V$ and that there
      is a point $x_{0}\in V$ whose $G^{c}$ orbit is dense. We also want
      to suppose  that the stabiliser
      $\Gamma\subset G^{c}$ is finite. Then the orbit is a copy of
      $G^{c}/\Gamma$ in $V$ and  the complement
      is an analytic subvariety  (which must contain a divisor if $X$
      is Kahler). Of course the case of a toric manifold fits into this picture,
      except that in that case we can assume $\Gamma$ is trivial (but see
      the further discussion below). In the next section we will study a
      particular example of this set-up: the Mukai-Umemura manifold.
      
      Now there is no loss of generality in supposing that $\Gamma$ lies
      in the compact group $G$ and we can study $G$-invariant Kahler metrics
      on $V$. Over the dense orbit these can be represented by Kahler potentials
      $\Phi$ on $G^{c}$ which are invariant under the two groups $G$ (acting
      by left multiplication) and $\Gamma$ (acting by right multiplication).
      In other words, $\Phi$ can be regarded as a function on the symmetric
      space $M=G^{c}/G$ which is invariant under the action of the finite  group
      $\Gamma$ on $M$. We will denote the corresponding function
      on $M$ by $\phi$.

      A finite group $\Gamma$ can enter in the toric case in slightly different
      way, but leading to the same conclusion. Suppose $\Gamma$ is a finite subgroup of $GL(n,\bZ)$ which preserves
      the polytope $P$ of a toric manifold $X$. (For example if $X$ is $\bC\bP^{n}$,
      so $P$ is the standard simplex, we can take $\Gamma$ to be the permutations
      of the $n$ coordinates.) 
      Then there is a group  $\hat{T}$  which fits into a split exact
sequence
\begin{equation} 1 \rightarrow T \rightarrow \hat{T} \rightarrow \Gamma \rightarrow
1 \end{equation} and which  acts on $X$. As a toric manifold, we know that
we can represent $T$-invariant Kahler metrics on $X$ by potentials $\phi$ on $\bR^{n}$, but now we can further restrict to $\hat{T}$-invariant metrics and
these correspond to $\Gamma$-invariant functions $\phi$, for the natural
action of $\Gamma$ on $\bR^{n}$ (of course, this copy of $\bR^{n}$ is really
the dual of that containing $P$).

\

    We now develop the local Kahler differential geometry in this
    situation, working in terms of a function $\phi$ on the symmetric space
    $M$.    This has a standard connection on its tangent bundle, which is
    the Levi-Civita connection for any $G^{c}$-invariant metric. Thus we
    have a Hessian operator $\nabla^{2}$ taking functions on $M$ to  sections
    of $s^{2}(T^{*}M)$.  The tangent space of $V$ at a point $g x_{0}$ can be
    identified with the complexification of the tangent space of $M$ at the
    point $Gg$. Thus we have an identification with the symmetric tensors
    $s^{2}(T^{*} M)$ at $Gg$ with a subspace of $\Lambda^{1,1} TG^{c}$ at $g$. This just corresponds to embedding the real symmetric matrices in the
complex Hermitian matrices. 
\begin{lem}
Under this identification for any function $\phi$ on $M$ and corresponding
function $\Phi$ on $G^{c}$ the form $i\partial \dbd \Phi$ corresponds to $\nabla^{2}\phi$.
\end{lem}
We can see this as follows. First note that in the toric case this is just
what we have seen when we identify the Kahler metric with the Hessian
$\phi^{ab}$.  For the general case, there is no loss in working
at the point $g=1$. To evaluate $\nabla^{2}\phi$ on a tangent
vector $v$  we take the geodesic $\gamma(t)$ in $M$ starting with
initial velocity $v$. Then 
$$   \nabla^{2}\phi(v)= \frac{d^{2}}{dt^{2}} \phi(\gamma), $$
evaluated at $0$. Now geodesics in $G^{c}/G$ through the identity coset correspond to $1$-parameter subgroups in $G^{c}$ so we have a homomorphism
$\tilde{\gamma}: \bC \rightarrow G^{c}$, such that $\gamma(t)= K\tilde{\gamma}(it)\in
M$. Then we are essentially reduced to the toric case, restricting to this
$1$-parameter subgroup.

Thus the local Kahler geometry in this situation reduces to the study of
{\it convex} functions on $M$ which, by definition, are those functions $\phi$
 with $\nabla^{2} \phi>0$ at each point. Equivalently, they are functions
 which are convex along geodesics in $M$. Of course this is a generalisation
 of the case when $M=\bR^{n}=
 T^{n}_{c}/ T^{n}$. We can go on to write out the equations we want to solve
 explicitly in this framework. The Kahler-Einstein equation, in the Fano
 case, is
 $$\det \nabla^{2} \phi= e^{-\phi}. $$
 For the scalar curvature; given a convex function $\phi$, we define an operator
 $$\Delta_{\phi}(f)= (\nabla^{2}\phi)^{-1}.\nabla^{2} f, $$
 where $(\nabla^{2}\phi)^{-1}$ is the quadratic form on $T^{*}M$ induced
 by the nondegenerate quadratic form $\nabla^{2}\phi$ on $TM$, in the usual
 way, and the dot denotes the contraction between $s^{2}TM$ and $s^{2} T^{*}M$.
 Then the scalar curvature of the Kahler metric defined by $\Phi$ is
 $$ S= \Delta_{\phi} (\log \det \nabla^{2} \phi). $$
 
 Notice that these local constructions make sense on any manifold equipped
 with a connection and volume form. 

There are some important differences between this theory in the case of
a semi-simple group $G$ and that in the abelian, toric, case.
\begin{itemize}
\item When we go beyond the local differential geometry we need to consider
a class of \lq\lq admissible'' functions $\phi$ which define metrics which
extend smoothly to $V$. This imposes some asymptotic growth conditions on
$\phi$ (as in the toric case) but these can be  more complicated, since they
encode the structure of the compactification.
\item In the toric case the local equations are affine invariant, but there
is no substitute for the affine group in the semi-simple case. In the semi-simple
case we have a preferred metric which changes the character of the theory.
\item The geometry of $M$ in the semi-simple case has negative curvature,
reflecting the non-abelian nature of $G$. This makes a radical difference
to arguments involving volumes of balls etc.
\end{itemize}

Again, there seems to the author to be a lot of scope for development of this theory.
For example one could consider a function $w$ on a Riemannian manifold of negative
curvature which satisfies a differential inequality
$$  \det \nabla^{2} w \geq e^{-w}, $$
and try to establish analogs of the results proved by Wang and Zhu in the toric case.

      \section{The Mukai-Umemura manifold and its deformations}
      
      The  first part of this section  gives an account, not aimed at algebraic
      geometry specialists, of a very interesting family of Fano $3$-folds,
      following Mukai.
      The basic references are \cite{kn:Mukai}, \cite{kn:MU}, but there are
      also many other relevant papers in the algebraic geometry literature.
      Then we go on to discuss the existence of Kahler-Einstein metrics
      on some manifolds in this family. 
      \subsection{Mukai's construction}
      
      We start with a $7$-dimensional complex vector space $V$ and write
      $Gr_{3}(V)$ for the Grassmann manifold of 3-dimensional subspaces of $V$.
      So $Gr_{3}(V)$ has dimension $3.(7-3)=12$. A form $\Omega\in \Lambda^{2}(V^{*})$
      defines a subset $Z_{\Omega}\subset Gr_{3}(V)$ consisting of the $3$-planes
      $P$ such that $\Omega\vert_{P}$ vanishes. In other language we consider
      the tautological rank 3 vector bundle $U\rightarrow Gr_{3}(V)$; the form $\Omega$
      defines a section $s_{\Omega}$ of $\Lambda^{2} U^{*}$ with zero set $Z_{\Omega}$.
      For generic $\Omega$ this zero set is a smooth subvariety of codimension
      $3$. Now let $\Omega_{1}, \Omega_{2}, \Omega_{3}$ be three such forms
      and consider 
      $$  X=Z_{\Omega_{1}} \cap Z_{\Omega_{2}} \cap Z_{\Omega_{3}} \subset
      Gr_{3}(V). $$
      Of course this only depends on the $3$-plane $\Pi$ in $\Lambda^{2}
      V^{*}$ spanned by the $\Omega_{i}$, so we may sometimes write $X_{\Pi}$.
      Obviously there is a Zariski-open subset ${\cal U}$ in the Grassmannian
      $Gr_{3}(\Lambda^{2} V^{*})$ of $3$-planes $\Pi$ such that $X_{\Pi}$ is a smooth subvariety
      of dimension $12-3.3=3$. This set ${\cal U}$ is non-empty, as we will  see
      later. The group $SL(V)$ acts on the whole construction and obviously
      different subspaces $\Pi$ which lie in the same $SL(V)$ orbit define
      isomorphic manifolds $X_{\Pi}$, so we get a {\it set} of equivalence classes
      of manifolds constructed in this way,  parametrised by the
      quotient
      ${\cal U}/SL(V)$. (Mukai shows further that this
      parametrisation is effective: i.e. $X_{\Pi_{1}}$ is isomorphic to $X_{\Pi_{2}}$
      if and only if $\Pi_{1},\Pi_{2}$ lie in the same $SL(V)$ orbit. Moreover,
      he shows that all \lq\lq prime Fano $3$-folds of genus 12'' arise in
      this way.)

      We compute the canonical bundle $K_{X}$ of the variety $X=X_{\Pi}$ for some
      $\Pi\in {\cal U}$. We have 
      $$   \Lambda^{2} U^{*}= U \otimes H $$
      where $H$ is the ample line bundle $\Lambda^{3} U^{*}$. So, writing
      $\det$ for the  the top exterior power of a vector bundle, we have
      $$  \det \Lambda^{2} U^{*}= H^{\otimes 2}.$$ The tangent
      bundle of the Grassmannian at a $3$-plane $P\subset V$ can be identified
      with $P^{*}\otimes V/P$. So 
      $$  \det TGr_{3}= H^{\otimes 7}. $$
      Now since the tangent bundle of $X$ is the kernel of a surjective map
      from $TGr_{3}(V)$ to $\Lambda^{2}U^{*}\oplus \Lambda^{2} U^{*} \oplus \Lambda^{2}
      U^{*}$ we have
      $$  K_{X}^{-1}= \det TX= H^{\otimes(7-3.2)}=H. $$
      Thus $X$ is a  Fano manifold. 
      The sections of $H$ over $Gr_{3}$ give the Plucker embedding 
      $$Gr_{3}(V)\rightarrow \bP(\Lambda^{3} V)= \bP^{34}$$
      For any $3$-form $A\in \Lambda^{3} V^{*}$ we get a hyperplane section
      $Y_{A}\subset Gr_{3}(V)$ which just consists of the $3$-planes $P$ such
      $A\vert_{P}=0$. By definition this occurs if $P$ is in $X$ and
      $A$ is in the image of the wedge product map $\Pi\otimes V^{*}\rightarrow
      \Lambda^{3} V^{*}$. We expect this map to have an image of dimension
      $7.3=21$ in which case the image of the composite
      $$  X\rightarrow Gr_{3}(V) \rightarrow \bP^{34}$$ lies in a linear subspace
      $\bP^{34-21}=\bP^{13}$. Certainly this map is defined by sections of
      $K_{X}^{-1}$, we will see later that $H^{0}(X,K_{X}^{-1})$ has dimension
      $14$ and that this embedding is  that given by the anticanonical
      system.
      
      To make this more concrete we show now that
       $X$ is a {\it rational} variety; that is, we construct an explicit
       parametrisation of a dense open set in $X$. Suppose we have a pair
       of $3$-dimensional subspaces $P_{0},Q_{0}\subset V$ with $P_{0}\cap Q_{0}=0$.       We ask what $3$-planes $P$ in the $6$-dimensional subspace $P_{0}\oplus
      Q_{0}$ lie in $X$. In matrix notation, we can write the restriction
      of a  form $\Omega$ to $P_{0}\oplus Q_{0}$ as
       $$ \left(\begin{array}{cc} \sigma &A\\ -A^{T}& \tau\end{array}\right)$$
        Now consider the $3$-dimensional subspaces $P$ which arise as the
        graphs of linear maps $M:P_{0}\rightarrow Q_{0}$. The condition becomes
        \begin{equation}\sigma + M^{T} \tau M + (AM- (AM)^{T}) =0.\end{equation}
         So our three
        forms $\Omega_{i}$ give us three triples $A_{i}, \sigma_{i}, \tau_{i}$
        and we have three equations of the form (33) to solve to find a point
        of $X$. We have $9$ unknowns: the entries of the matrix $M$.  The left hand side of (33) takes values in the $3$-dimensional
        space of skew symmetric $3\times 3$ matrices so we obtain a total of
        $3.3=9$ equations in these $9$ unknowns and we expect a finite number of solutions. These equations are
        quadratic and one can  solve them explicitly,  to see that there are generically two solutions.
        However it is easier to suppose  that we are in the case when $P_{0}$ itself lies
        in
        $X$. This means that all the $\tau_{i}$ are zero, so the equations
        (33) become linear. Generically this system of $9$ linear equations
        in $9$ unknowns is nondegerate and there is a unique solution.
        Now suppose we have found one point $P_{0}$ in $X$ and consider the
        space of $6$-planes in $V$ which contain $P_{0}$. This is a copy
        of projective $3$-space $\bP^{3}$. Given a point in $\bP^{3}$, that
        is to say a 6 dimensional subspace $E$ of $V$, we choose a complementary
        subspace to write is as $E= P_{0}\oplus Q_{0}$. Then we can proceed
        as above and, by solving linear equations, find the points of $X\cap
        Gr_{3}(E)$. Generically there is just one, $P_{E}$ say,  different from the original
        $P_{0}$. Conversely for any $P'\in X$ the sum $P\oplus P'$ lies
        in a $6$-dimensional subspace. Of course there will be various exceptional
        cases, but the upshot is that we get a birational map from $\bP^{3}$
        to $X$ which takes a subspace $E$ containing $ P_{0}$ to $P_{E}$.
        
      \

We now consider a special manifold in this family. Take the vector space $V$ to be the
sixth symmetric power $s^{6}$ of the fundamental representation of $SL(2,\bC)$.
Then $\Lambda^{2} V^{*}= \Lambda^{2} s^{6}$  decomposes into distinct irreducible representations
$$  \Lambda^{2} s^{6}= s^{10}\oplus s^{6}\oplus s^{2}. $$
The $s^{2}$ summand is a $3$-plane $\Pi_{0}$ invariant under $SL(2,\bC)$, so
there is a natural $SL(2,\bC)$ action on the corresponding variety, the {\it
Mukai-Umemura manifold}, $X_{0}=X_{\Pi_{0}}$. We will see below that
$X_{0}$ admits a Kahler-Einstein metric.
The representation $s^{6}$ has a standard invariant symmetric form $(\ ,\
)$ and the inclusion $s^{2}\rightarrow \Lambda^{2} s^{6}$ is just the map
from the Lie algebra of $SL(2,\bC)$ given by the action on $s^{6}$. This
comes down to saying that a $3$-plane $P$ is in $X_{0}$ if and only if
\begin{equation} (\delta p, q)=0 \end{equation}
for all $p,q\in P$ and $\delta\in \mathfrak{s}\mathfrak{l}_{2}$. 
Notice that the action of $SL(2,\bC)$ on all the spaces involved actually
factors through $PSL(2,\bC)$.

Identify the projectivisation of the fundamental representation $s^{1}= \bC^{2}$
with the standard round sphere and fix an {\it icosahedron}, which can be
regarded as a set of 12 vertices in  in this sphere.
Thus we get a symmetry group $\Gamma\subset SO(3)\subset PSL(2,\bC)$ of order $60$. There is a simple way to see that the $7$-dimensional representation $s^{6}$ of $PSL(2,\bC)$
becomes reducible when restricted to $\Gamma$. There are $6$ pairs of
antipodal vertices and for each such pair $p,\overline{p}$ we have a $1$-dimensional
subspace consisting of polynomials which vanish to order $3$ at $p,\overline{p}$.
The sum of these $6$ subspaces is obviously invariant under $\Gamma$ and
is a proper subspace of $s^{6}$ since it has codimension at least $1$. A
little calculation shows that this invariant subspace is of dimension $3$
and satisfies the criterion (34). So this subspace gives a point $P_{0}$ in
$X_{0}$ fixed by $\Gamma$. On the other hand the stabiliser of $P_{0}$ is
obviously not the whole of $SO(3)$ and, since there is no finite subgroup
of $SO(3)$ strictly larger than $\Gamma$, the stabiliser must be exactly $\Gamma$.

Now go back to the wedge product $P_{0}\wedge  V^{*}\rightarrow
\Lambda^{3} V^{*}$. In terms of representations this is an $SL(2,\bC)$-map
$$    s^{2}\otimes s^{6}\rightarrow \Lambda^{3} s^{6}. $$
It is an exercise in representation theory to show that
$$ \Lambda^{3} s^{6} = s^{12}\oplus s^{8}\oplus s^{6}\oplus s^{4}\oplus s^{2}\oplus
s^{0}\oplus s^{0}. $$
So comparing with 
$$   s^{2}\otimes s^{6}= s^{8}\oplus s^{6}\oplus s^{4}\oplus s^{2}\oplus
s^{0}$$
we see that the embedding $X_{0}\subset Gr_{3}(V)\subset \bP(\Lambda^{3} V)$
gives rise to an $SL(2,\bC)$-equivariant embedding
  \begin{equation}  X_{0} \rightarrow \bP(s^{0} \oplus s^{12}). \end{equation}
  In other words, by our identification of the anticanonical bundle $K^{-1}$ we have
  $$  H^{0} (X_{0}, K^{-1}) = s^{0} \oplus s^{12}, $$
  as a representation of $SL(2,\bC)$.In  particular, there is an $SL(2,\bC)$-invariant section $\sigma$ of $K^{-1}$. Explicitly, if we identify $\Lambda^{3}
  s^{6}$ with $\Lambda^{4} s^{6}$ then $\sigma$ corresponds to the $4$-form
  on $V=s^{6}$ defined as follows. We choose any orthonormal basis $\Omega_{1},
  \Omega_{2}, \Omega_{3}$ of $P_{0}$ and write down the $4$-form
  $$ *\sigma= \Omega_{1}^{2}+ \Omega_{2}^{2} + \Omega_{3}^{2}. $$
  
  \
  
  In this way, we  get another description of the manifold $X_{0}$. Our point $P_{0}\in
  X_{0}$ cannot lie in the zero set of $\sigma$ (since its orbit is $3$-dimensional).
  So, in the embedding (35), we have
  $$  P_{0}= [ 1, v_{0}] \in \bP(\bC\oplus s^{12}). $$
  Thus $v_{0}$ is an element of $s^{12}$ whose stabiliser in $ PSL(2,\bC)$
  is exactly $\Gamma$.Now there is an obvious element of the projective space
  $\bP(s^{12})$ with stabiliser $\Gamma$, just the configuration of vertices
  of the icosahedron, regarded
  as an element of the symmetric product. Since $\Gamma$ is a perfect group
  it must act trivially on the  corresponding line in $s^{12}$, so we get
  a vector in $s^{12}$  with stabiliser $\Gamma$. It is easy to see that, up
  to a multiple, this in the only element of $s^{12}$ with stabiliser $\Gamma$,
  and thus we have  identified $v_{0}$. Then we can simply define $X_{0}$ to be the
  closure in $\bP(\bC\oplus s^{12})$ of the $PSL(2,\bC)$-orbit of $v_{0}$
  in $s^{12}$. (Here we are regarding the vector space $s^{12}$ as being
  a subset of the projective space $\bP(\bC\oplus s^{12})$ in the familiar way.)

  In this description, the intersection of $X_{0}$ with the hyperplane at
  infinity
  $$ D= \bP (s^{12}) \subset \bP(\bC\oplus s^{12}), $$
  is, by definition, the zero set of the invariant section $\sigma$ of $K^{-1}$.
  Consider a $1$-parameter subgroup $\lambda_{t}$  in $PSL(2,\bC)$. Thus
  we have a pair of distinct point $z_{+}, z_{-}$ such that when $t$
  is large positive  the map $\lambda_{t}$ contracts most of the sphere to
  a small neighbourhood of $z_{+}$, and when $t$ is large negative to a small
  neighbourhood of $z_{-}$. If $y_{1}, \dots y_{12}$ is any configuration
  of distinct points it is not hard to see that the limit as $t\rightarrow
  \infty$ of
  $$ \lambda_{t}(\underline{y})= \left(\lambda_{t}(y_{1}), \lambda_{t}(y_{2} \dots
  \lambda_{t}(y_{12})\right)$$
  in the symmetric product $\bP(s^{12})$ is either $12 z_{+}=(z_{+}, z_{+}, \dots,z_{+})$
  (in the generic case) or $11 z_{+}+ z_{-}=(z_{+}, \dots, z_{+}, z_{-})$ (in the case when
  one of the $y_{i}$ is $z_{-}$). Using this, Mukai and Umemura show that
  the divisor at infinity $D$  consists precisely of the union of points
  of the form $12 z_{+}$ or $11 z_{+}+ z_{-}$ in $\bP(s^{12})$. It is easy
  to identify
  this geometrically. The points of the form $12 z_{+}$ make up the {\it
  rational normal curve} in $\bP(s^{12})$. Our divisor $D$ is the surface
  swept out by the lines in $\bP(s^{12})$ tangent to the rational normal
  curve. As a set we can identify $D$ with $\bP^{1} \times \bP^{1}$: we just
  map $(z_{+}, z_{-}) \in \bP^{1}\times \bP^{1}$ to $11 z_{+}+ z_{-} \in
  D$. But the surface $D$ is singular and a more precise statement is that
  the map above is a holomorphic map $\nu: \bP^{1}\times \bP^{1} \rightarrow D$ which
  is the normalisation of $D$. The singular set of $D$ is the image of the
  diagonal in $\bP^{1} \times \bP^{1}$, and it is easy to check that the singularity
  has the form of a  cusp transverse to the diagonal. That is to say,
  we can choose  local co-ordinates $z_{1} z_{2}, z_{3}$ in $X_{0}$ around
  a singular point of $D$ such that $D$ is defined by the equation $z_{1}^{2}=
  z_{2}^{3}$. 
  
  We now have a rather explicit description of $X_{0}$, as the compactification
  of $PSL(2,\bC)/\Gamma$ formed by adjoining the divisor $D$. We can use
  this to compute the action of $PSL(2,\bC)$ on all of the spaces of sections
  $H^{0}(X_{0}, K^{-p})$. For the pull back $\nu^{*}(K^{-1})$ is isomorphic
  to the line bundle ${\cal O}(11,1)$ over $\bP^{1}\times \bP^{1}$. We can
  regard the structure sheaf of $D$ as a subsheaf of  that of $\bP^{1}\times
  \bP^{1}$. From the local model of the singularity along the diagonal one
  sees that the quotient can be identified with sections of ${\cal O}(2)$
  along the diagonal. This means that $H^{0}(D, K^{-p}\vert_{D})$ is the
 kernel of a map $H^{0}(\bP^{1}\times \bP^{1}; {\cal O}(11p,p))\rightarrow
 H^{0}(\bP^{1}; {\cal O}(12p-2))$. As representations of $PSL(2,\bC)$ this
 is a map
 $$  s^{11p}\otimes s^{p}\rightarrow s^{12p-2}. $$
 Now $$s^{11p}\otimes s^{p}= s^{12p}+ s^{12p-2} \dots \oplus s^{10p}$$
 and the map above is just the projection to the second factor. So
 $$  H^{0}(D; K^{-p}\vert_{D})= s^{12p}\oplus s^{12p-4}\oplus s^{12p-6}\dots
 \oplus s^{10p+2}\oplus s^{10p}. $$
 Then the exact cohomology sequence of 
 $$  0\rightarrow K^{-(p-1)}\rightarrow K^{-p} \rightarrow K^{-p}\vert_{D}\rightarrow
 0$$
 together with Kodaira vanishing on $X_{0}$ gives
 $$  H^{0}(X_{0}, K^{-p})= H^{0}(X_{0}, K^{-(p-1)})\oplus s^{12p}\oplus s^{12p-4}\dots
 s^{10p}, $$
 and inductively we get a description of each $H^{0}(X_{0}, K^{-p})$. Thus
 $$  H^{0}(X_{0}, K^{-1}) = s^{0} \oplus s^{12}, $$
 $$  H^{0}(X_{0}, K^{-2})= s^{0}\oplus s^{12} \oplus s^{24}\oplus s^{20}.
 $$
 
 For $p\geq 6$ we get multiplicities: $H^{0}(X_{0}, K^{-6})$ contains two
 copies of $s^{60}$. This illustrates the difference with the multiplicity-free case discussed
 above. (Although since the multiplicities are small  until
 $p$ becomes quite large, once is tempted to think of $X_{0}$ as  being \lq\lq
 close'' to multiplicity-free. )

      \subsection{Topological and symplectic picture}
      
      We will now get another  explicit picture of $X_{0}$, taking
      the point of view of symplectic geometry. Recall that all Kahler metrics in
      the cohomology class $c_{1}(X_{0})$ define equivalent symplectic structures,
      so we have a well-defined symplectic manifold $(X_{0}, \omega)$ with an $SO(3)$-action. Thus we have an equivariant moment map
$$  \mu: X_{0} \rightarrow \bR^{3}= \Lie(SO(3))^{*}. $$
whose image is clearly a ball in $\bR^{3}$. We can understand the structure
of this moment map by restricting to a subgroup $S^{1}\subset SO(3)$, say
that corresponding to the $x_{1}$-axis in $\bR^{3}$. Then the Hamiltonian
$H$
for this circle action on $X_{0}$ is the composite of $\mu$ with projection
to the $x_{1}$-axis.  The critical points of $H$ are the fixed points of the
circle action and we can find these explicitly.  We can suppose that our
circle subgroups corresponds to the standard action of
 $$\left(\begin{array}{cc}\lambda^{1/2}&0\\0&\lambda^{-1/2}\end{array}\right),$$
acting on $\bC^{7}$ with weights $\lambda^{3}, \dots \lambda^{-3}$. We write
$e_{i}$ for the basis vector belonging to the weight $\lambda^{i}$.
This induces an action on the Grassmannian $Gr_{3}(V)$ whose fixed points are just
invariant $3$-dimensional subspaces of $\bC^{7}$ and these are just the spans
$P_{ijk}=\langle e_{i}, e_{j}, e_{k}\rangle$
for distinct $i,j,k$. By checking the 35 different cases, or otherwise, one
finds that the only $P_{ijk}$ which satisfy the criterion (34) to lie in $X_{0}$
are $P_{123}, P_{023}, P_{0-2-3}, P_{-1-2-3}$. There is an action of the
Weyl group $\{\pm1\}$ on the whole situation which commutes, up to sign, with the circle
action, takes $H$ to $-H$  and takes $P_{ijk}$ to $P_{-i-j-k}$. So there
are four fixed points of the circle action but to analyse the local structure
around them
it suffices to consider the two cases $P_{123}, P_{023}$. Notice that, by
considering $H$ as a Morse function we immediately see that $X_{0}$ has the
same additive homology as $\bC\bP^{3}$. Notice also that the value of $H$
at a critical point is just given by the weight of the action on the fibre
of $K^{-1}$ over this point, which is just $i+j+k$ at $P_{ijk}$.

We next compute the weights of the circle action on the tangent spaces at the fixed points. This
is similar to the calculation of the canonical bundle. At a fixed point the
tangent space $TX_{0}$, viewed as a representation of $S^{1}$, can be written
as the formal difference 
$$  TGr_{3}(V)-\left( \Lambda^{2} U^{*} \otimes\Lie(SO(3)) \right). $$
Computing the weights of these two terms and subtracting we find that the
weights of the action on the tangent space at $P_{123}$ are $(1,2,3)$ and
on the tangent space at $P_{023}$ are $(1,-1,5)$. In either case the orbit
of the fixed point is a copy of $SO(3)/S^{1}=S^{2}$ in $X_{0}$ and the weight
$1$ in the action on $TX_{0}$ just corresponds to the tangent space of this
orbit. The weights normal to the orbit are $(2,3)$ in the case of $P_{123}$
and $(-1,5)$  in the case of $P_{023}$.

With these calculations we can get a good picture of the map $\mu$. Write
$\Sigma, \Sigma'$ for the orbits of $P_{123}$ and $ P_{023}$ respectively. Then $\mu$ restricts to an $SO(3)$-equivariant equivalence
between $\Sigma$ and the sphere of radius $1+2+3=6$ in $\bR^{3}$ and between
$\Sigma'$ and the sphere of radius $0+2+3=5$. The image of $\mu$ is the ball
of radius $6$ and the critical values of $\mu$ are precisely these two spheres.
So $\mu$ is a fibration away from these spheres. For $\ux\in \bR^{3}$, write $F_{\ux}$ for the
preimage $\mu^{-1}(\ux)$. If $\vert \ux\vert\neq 5,6$ the fibre $F_{\ux}$
is a $3$-manifold. If also $\vert \ux\vert>0$ then this $3$-manifold has
a natural circle action defined by the circle subgroup of
$ SO(3)$ fixing $\ux$. When $\ux=0$ the fibre has an $SO(3)$ action. As $\ux$ varies in $\bR^{3}$ the fibre only \lq\lq
changes''---in the obvious sense---when $\vert \ux\vert$ crosses the special
values $5,6$. Thus we understand the full topological picture if we understand
the  changes in the fibre as $\ux$ moves along the positive $x_{1}$-axis,
say. Let $V\subset X_{0}$ be the pre-image by $\mu$ of the positive $x_{1}$-axis. This
is a smooth $4$-manifold, with a circle action, and the fibres $F_{\ux}$, for
$\ux$ on the axis, are the level sets of the Hamiltonian $H$, restricted
to $V$. Then we have the usual Morse-theory description of these changes,
from the Hessian of $H$ on $V$, which is determined by the weights of the
circle action. As $\ux$ moves across the point $(6,0,0)$ the situation is
modelled by the level sets
   $$  2 \vert z_{1}\vert^{2} + 3 \vert z_{3}\vert^{2} = \epsilon, $$
   for $(z_{1}, z_{2})\in \bC^{2}$, with the circle action of weight $(2,3)$.
   Thus the fibre changes from the empty set to a $3$-sphere with an action
   given by these weights. As $\ux$ moves across the point $(5,0,0)$ the
   situation is modelled, locally,  by the level sets
   $$ -\vert z_{1}\vert^{2} + 5\vert z_{2} \vert^{2} =\epsilon, $$
   with the circle action of weight $(-1,5)$.
The effect on the fibres is
to perform a \lq\lq Dehn surgery'' on an $S^{1}$-orbit. Thus the fibres $F_{\ux}$
for $\vert \ux\vert<5$
are obtained by performing this surgery on a knot $\Gamma\subset S^{3}$.
Now $\Gamma$ is a free orbit of the $(2,3)$ action so it is the $(2,3)$ \lq\lq
torus knot'' which is just a trefoil. To nail down the Dehn surgery completely we need to specify a framing
of the knot but this is determined by the fact that the linking number of
a nearby orbit with $\Gamma$ is the weight $5$, from which one concludes that the
framing is $+1$. This is a well-known description of the Poincar\'e homology
sphere (the result of $+1$-surgery on a trefoil), and ties in with our previous
discussion since the fibre $F_{0}$ is the $SO(3)$-orbit  $SO(3)/\Gamma$.
(Another way of expressing this is that the fibres $F_{\ux}$ are Seifert-fibred
$3$-manifolds: for $5<\vert \ux\vert<6$ we have two multiple fibres with
multiplicity $(2,3)$ and the surgery across $\vert \ux\vert=5$ introduces another
multiple fibre with multiplicity $5$, so for $\vert \ux\vert<5$ we get the
Seifert manifold with multiplicities $(2,3,5)$, which is another well-known
description of the Poincar\'e manifold.)

\

It is interesting to match this picture up with the algebro-geometric description.
This illustrates the general theory of Kirwan \cite{kn:Ki}. The $2$-sphere $\Sigma$
at which $\vert \mu\vert$ attains its maximal value $6$ is a holomorphic
sphere in $X_{0}$: it is just the rational normal curve in our divisor $D\subset
\bP(s^{12})$. The other sphere $\Sigma'$ is not holomorphic. It is a critical
manifold for the function $\vert \mu\vert^{2}$ on $X_{0}$ and the divisor $D$ appears
as the associated \lq\lq ascending set'': the closure of the set of points
which flow to $\Sigma'$ under the decreasing gradient flow of $\vert \mu\vert^{2}$.
In our description of $D$ as $S^{2}\times S^{2}$ the holomorphic curve $\Sigma$
is the diagonal and $\Sigma'$ is the \lq\lq anti-diagonal''consisting of pairs
of antipodal points. One can also see the cusp singularity in $D$, transverse
to $\Sigma$, from the weights $(2,3)$ of the circle action on the normal bundle.

Notice that if we write $\bC\bP^{3}=\bP(s^{3})$, for the $4$-dimensional
representation $s^{3}$ of $SU(2)$, the moment map $\mu:\bC\bP^{3}\rightarrow
\bR^{3}$ for the action gives a description of $\bC\bP^{3}$ very similar
to that above. In this case $\mu^{-1}(0)$ is $SO(3)/H$ where $H\subset SO(3)$\ is the group
of  symmetries of an equilateral triangle, and we see this $3$-manifold described as the
Seifert fibration with multiple fibres $(2,2,3)$.

      \subsection{Deformations}

      Here we study the deformations of Mukai's construction about the special
      solution $X_{0}$. Recall that a manifold in this family is specified
      by a $3$-plane in $\Lambda^{2} \bC^{7}$. We start with the $3$-plane
      $s^{2} \subset \Lambda^{2} s^{6}= s^{10}\oplus s^{6}\oplus s^{2}$. The tangent space of the Grassmannian
      at this point is given by the linear maps from $s^{2}$ to the complementary
      subspace $s^{10}\oplus s^{6}$, that is (using the fact that all these
      representations are isomorphic to their duals)
         $$TGr_{3}(\Lambda^{2} \bC^{7})= (s^{10}\oplus s^{6})\otimes s^{2}=
         s^{12}\oplus 2 s^{8} \oplus s^{6}\oplus s^{4}. $$
         The action of the group $SL(\bC^{7})=SL(s^{6})$ gives a linear map
         $$  \mathfrak{s}\mathfrak{l}(7) \rightarrow TGr_{3}(\Lambda^{2}\bC^{7}), $$
         which we know has kernel the Lie algebra $\mathfrak{s}\mathfrak{l}(2)$ of the stabiliser.
         Now the Lie algebra of $GL(\bC^{7})=GL(s^{6})$ is $$s^{6}\otimes s^{6}=s^{12}\oplus
         s^{10}\oplus s^{8}\oplus s^{6}\oplus s^{4}\oplus s^{2}\oplus s^{0}$$
         so the Lie algebra of $SL(s^{6})$ is $s^{12}\oplus \dots s^{2}$.
         It is clear then that the quotient of the tangent space by the tangent
         space to the orbit is just $s^{8}$, as a representation of $PSL(2,\bC)$.
         By general theory there is an equivariant slice: a  $PSL(2,\bC)$                 equivariant embedding $j$ from
         a neighbourhood of $0$ in $s^{8}$ into $Gr_{3}(\Lambda^{2})$, mapping
         $0$ to our fixed subspace $s^{2}$, such that two points $j(p), j(q)$
         in  the same $SL(7)$ orbit if and only $p,q$ are in
         the same $PSL(2)$ orbit. In fact, although we do not really need
         this,  what we are describing is the versal
         deformation of $X_{0}$, so $H^{1}(TX_{0})= s^{8}$, as a representation
         of $PSL(2,\bC)$.

         One can gain a lot of insight from this simple calculation. The
         structure of the orbits of $PSL(2,\bC)$ on $s^{8}$ (or any $s^{p}$)
         is a standard example in Geometric Invariant Theory.
         There are five cases
         \begin{enumerate}
         \item The trivial orbit $\{0\}$.
         \item
         The orbits of  polynomials having no zero of multiplicity $\geq 4$. These are closed in $s^{8}$.
\item The orbit of polynomials having two distinct zeros, each of multiplicity
four. This orbit is closed and each point in it has stabiliser $\bC^{*}\subset
PSL(2,\bC)$.
\item The orbits of polynomials having a zero of multiplicity four and other
zeros each of multiplicity less than four. These orbits are not closed but there closure contains
the orbit of type (3).
\item The orbits of polynomials having a zero of multiplicity $\geq 5$. these
are not closed and contain $0$ in their closure.
\end{enumerate}

This is the source of the famous example of Tian of Fano manifolds without
Kahler-Einstein (or Ricci soliton) metrics \cite{kn:T2}. Tian shows that the manifolds
corresponding to any  $PSL(2,\bC)$ orbit of type
(5) cannot have such metrics. Tian's general results also show the same for the manifolds
corresponding to orbits of type (4). Tian's results are of course deep and
difficult but we note now  that a weaker statement is rather obviously true.
For this we need to recall some background.

In  general, the linearisation of the Kahler-Einstein equations on a complex
manifold $Z$ at a solution
$\omega_{0}$ is given by the self-adjoint operator $\Delta+1$ and, much as
we have seen in Section 3, the kernel of this can be identified with the
Lie algebra of the isometry group $G$ of $\omega_{0}$. Suppose we have
a $G$-equivariant deformation of $Z_{0}$: i.e. a complex manifold ${\cal
Z}$ with a $G$-action, an action of $G$ on a ball $B\subset \bC^{m}$ and
a $G$-equivariant submersion $\pi:{\cal Z}\rightarrow B $. In this situation
we automatically get  \lq\lq local actions'' of the complexified group $G^{c}$
on ${\cal Z}$ and $B$, compatible with $\pi$.
The standard \lq\lq
Kuranishi method'', which depends only on the formal properties of the
situation, yields the following structure (after possibly restricting to
a smaller ball $B$).
\begin{itemize}\item A $G$-invariant family of Kahler metrics
$\omega_{t}$ on the fibres $Z_{t}=\pi^{-1}(t)$ such that  $\omega_{t}$ is
isometric to $\omega_{t'}$ if and only if $t$ and $t'$ are in the same $G$-orbit.
\item A smooth 
map 
$\nu:B\rightarrow \mathfrak{g}^{*}$, equivariant for the action of $G$ on
$B$ and the co-adjoint action on $\mathfrak{g}^{*}$, such that $\omega_{t}$
is Kahler-Einstein if and only if $\nu(t)=0$.

\end{itemize}

Now in this general situation we can see that, if the $G$-action on $B$ is
non-trivial the map $\nu$ cannot be identically zero. For if $t,t'$
are in the same  orbit of the local $G^{c}$ action on $B$ then $Z_{t}$ and $Z_{t'}$ are isomorphic complex
manifolds. But if $\nu(t)$ and $\nu(t')$ both vanish then $\omega_{t}$ and
$\omega_{t'}$ are Kahler-Einstein and, by the uniqueness of the Kahler-Einstein
solution, they must be isometric and this only happens if $t,t'$ are in the
same $G$-orbit. Thus what we see from this elementary argument is that as we deform $Z_{0}$ in the smooth family $Z_{t}$ we
cannot deform the metric $\omega_{0}$ in a smooth family of Kahler-Einstein
metrics, for all small $t$. Tian's much stronger result is that if the Futaki
invariant of $Z_{0}$ vanishes (say), and if  $0$ lies
in the closure
of the the $G^{c}$-orbit of a point $t\in B$ then $Z_{t}$ does not admit
any Kahler-Einstein metric at all. This is an example of the \lq\lq jumping
of structures'' phenomenon discussed
in Section 1: there are arbitrarily small deformations of $Z_{0}$ which are
equivalent to a different structure $Z_{t}$.

         Returning to our special case of the Mukai-Umemura manifold, we can see conversely  that there are  some deformations
         of $X_{0}$ which {\it do} admit Kahler-Einstein metrics. 
         The general theory of these \lq\lq obstruction maps''  $\nu$ is
         being developed by T. Br\"onnle, in his Ph.D thesis, but in this
         special case we can make some  simple deductions from symmetry arguments.
           Let $p$
         be a point in $s^{8}$ which is fixed by a subgroup $J\subset SO(3)$.
         Then $J$ acts on $\bR^{3}=\mathfrak{s}\mathfrak{u}(2)$ and if $\nu$
         is any equivariant map from $s^{8}$ to $\bR^{3}$ then $J$ must fix
         $\nu(p)$. So if the origin is the only point in $\bR^{3}$ fixed by $J$ then we must have $\nu(p)=0$. Consider, for example,
$$   p= C (z^{4}-\alpha w^{4})(w^{4}-\alpha z^{4}), $$
with any $\alpha, C \in \bC$. This is fixed by a dihedral group $J$ of order $8$ which has the desired property, so we see that the deformations corresponding
such elements of $s^{8}$ admit Kahler-Einstein metrics, for small $C$. For $\alpha,C\neq 0$
the element $p$ has a discrete stabiliser in $SO(3)$ and it follows that
the corresponding metrics have discrete isometry groups. But then the deformation
theory implies that {\it all} small deformations of  these manifolds admit
Kahler-Einstein metrics. So we conclude that there is a non-empty open set
in ${\cal U}$ where the manifolds admit Kahler-Einstein metrics.

      Taking $\alpha=0$ above we get a special family of deformations, admitting
      Kahler-Einstein metrics, where we can take $J=O(2)\subset SO(3)$. It
      follows that the corresponding manifolds have a $\bC^{*}$-action.
      We can see this family of manifolds explicitly
      as follows. Fix the action on $\bC^{7}$ with weights $\lambda^{3},\dots,\lambda^{-3}$
      as usual. Then we want to look at $3$-dimensional subspaces $\Pi$ of $\Lambda^{2}
      \bC^{7}$ preserved by the action and we just consider those on which
      the action has weights $1,0-1$. Now the weight $1$-subspace of $\lambda^{2}$
      has a basis $e_{3}\wedge e_{-2}, e_{2}\wedge e_{-1}, e_{1}\wedge e_{0}$ and our space
      $\Pi$ must contain a vector
      $$ u=  u_{3,-2}e_{3}\wedge e_{-2}+ u_{2, -1} e_{2}\wedge e_{-1}+ u_{1, 0} e_{1}\wedge e_{0},
      $$
      for scalars $u_{3, -2}$ etc.
      Similarly $\Pi$ must contain a vector
      $$ v=  v_{1,-1}e_{1}\wedge e_{-1} + v_{2, -2} e_{2}\wedge e_{-2} + v_{3, -3} e_{3}\wedge e_{-3}
      $$
      and a vector
      $$ w= w_{-3, 2} e_{-3}\wedge e_{2}+ w_{-2, 1} e_{-2}\wedge e_{1} + w_{-1, 0} e_{-1}\wedge
      e_{0}. $$
      The vector space $\Pi$ is determined by these three vectors $u,v,w$. The coefficients
      are not unique. We could change $u,v,w$ to $\mu_{1} u, \mu_{2}v, \mu_{3}
      w$. Also we could change our basis vectors $e_{i}$ to $\lambda_{i}
      e_{i}$ to give an equivalent $3$-plane. This would change the coefficients,
      for example $u_{3, -2}$ would change to $\lambda_{3}\lambda_{-2} u_{3,
      -2}$. However the expression
      $$ \tau =\frac{  u_{3, -2} w_{-3, 2} v_{1, -1}}{u_{2, -1} w_{-2, 1} v_{3, -3}}$$
      is invariant under all these changes and gives a \lq\lq modulus'' for
      this family. The Mukai-Umemura manifold has $\tau=1$. When $\tau$ is
      close to $1$ we have seen that the corresponding manifold admits a Kahler-Einstein
      metric. It seems likely that this true for all $\tau$ but, as far the
      author is aware, this is not known. It seems an interesting test
      case for future developments in the existence theory.

      \subsection{The $\alpha$-invariant}
      
      In this subsection we establish the fact used above, that the Mukai-Umemura manifold has a Kahler-Einstein
      metric\footnote{This material appeared in the preprint {\it A note
      on the $\alpha$-invariant of the Mukai-Umemura 3-fold} arxiv DG 07114357.},
      which is . For this we appeal to the theory of the $\alpha$-invariant,
      developed by Tian \cite{kn:T1}. 
     We begin by recalling the definition.
Let $Z$ be a Fano manifold on which a compact group $G$ acts by holomorphic
automorphisms and fix a $G$-invariant Kahler metric $\omega_{0}$ in the
cohomology class $-c_{1}(K_{Z})$. Let $\cP$ be the set of $G$-invariant Kahler
potentials $\psi$ on $X$ such that $\omega_{\psi}= \omega_{0}+i \partial
\overline{\partial} \psi>0$ and $\max_{Z} \psi=0$. Thus $\cP$ can be identified
with
the set of all $G$-invariant Kahler metrics in the given Kahler class. Let
$A\subset \bR$ be the set defined by the condition that $\beta\in A$ if there
exists a $C_{\beta}\in \bR$ such that
$$ \int_{Z} e^{-\beta \psi} d\mu_{0} \leq C_{\beta}, $$
for all $\psi\in \cP$. Here $d\mu_{0}$ is the volume form defined by the fixed
metric $\omega_{0}$. Then Tian sets
$$  \alpha_{G}(Z) = \sup\{\beta: \beta \in A\}, $$
and shows that this does not depend on the choice of $\omega_{0}$.
 He shows that $\alpha_{G}(Z) $ is always strictly positive and that if $\alpha_{G}(Z)>\frac{n}{n+1}$ then $Z$ has a Kahler-Einstein
metric. What we really show in this subsection is that if we take the Mukai-Umemura
manifold $X$ with the action of $SO(3)$ then,

\begin{thm}
The $\alpha$-invariant $\alpha_{SO(3)}(X_{0})$ is $5/6$.
\end{thm}

So, since $5/6>3/4$, Tian's theory proves the existence of a Kahler-Einstein
metric. 
We should say straightaway that this is not really a new result. Alessio
Corti has explained to the author that, given the facts above, it can be obtained from the more general
theories of \cite{kn:DKol}. But our argument is extremely simple and fits
well into the general framework of this article.
\

\


We will only write down the proof that $\alpha\geq 5/6$, which is what is
relevant to Corollary 1. The proof that $\alpha=5/6$ is an easy extension
of this.

\

\begin{lem}
There is an $M\in \bR$ such that
$$\int_{Z} \psi \ d\mu_{0} \geq -M$$
for all $\psi\in \cP$.
\end{lem}
This is a step in Tian's proof that $\alpha>0$ and we repeat his argument.
If $\psi\in \cP$ we have
$$  \Delta_{0} \psi = 2\Lambda( i\partial\overline{\partial} \psi)\geq
-2n. $$
Let $K$ be the Green's function for $\Delta_{0}$, so that for all functions
$f$ on $Z$
$$  f(x) = -\int_{Z} K(x,y) (\Delta_{0} f)(y) d\mu_{0}(y)+\frac{1}{V} \int
f(y) d\mu_{0}(y), $$
where $V$ is the volume of the manifold. With our sign conventions, $K$ is bounded below and, since we can change $K$ by the addition of a constant without
affecting the identity, we may suppose that $K\geq 0$. While $K$ is singular
along the diagonal it is  integrable in each variable. Let $x$ be the point
where $\psi$ vanishes. Then applying the Green's identity to $\psi$ we have
$$\int_{Z} \psi(y) d\mu_{0}(y) = V \int_{Z} K(x,y) \Delta_{0}\psi d\mu_{0}(y)\geq
-2n V \int_{Z} K(x,y) d\mu_{0}(y). $$
So we can take $$M= 2n V \max_{x}\int_{Z} K(x,y) d\mu_{0}(y). $$

\

\

For the rest of this section we work with  the Mukai-Umemura manifold, which
we denote by
$X$. Let $\sigma$ be the $SO(3)$-invariant  section of the anticanonical bundle $K^{-1}$ cutting
out the divisor $D$. There is a Hermitian metric  on this line bundle such
that the curvature of the associated unitary connection is $- i \omega_{0}$.
Set $$f_{0}=- \log \left(\vert \sigma \vert^{2}\right).$$
This is a smooth function on $X\setminus D$ and $i\partial \overline{\partial}
f_{0} = \omega_{0}$.
 
\begin{lem}
For any $\beta<\frac{5}{6}$ the function $\exp(\beta f_{0})$ is integrable.
\end{lem}

This is also standard. The integral in question is
$$  \int_{Z} \vert \sigma \vert^{-2\beta} d\mu_{0}. $$
By what we know about the singlarities of $D$, we can reduce to considering the
integrals
$$  \int_{B} \vert z^{2} - w^{3} \vert^{-2\beta} , $$
where $B$ is the unit ball in $\bC^{2}$ and $z,w$ are complex co-ordinates.
Let $T$ be the linear map $T(z,w)= (z/8, w/4)$ and for $r\geq 1$ set
$$ \Omega_{r}= T^{r}(B)\setminus T^{r-1}(B). $$
 Set
$$I_{r} = \int_{\Omega_{r}}  \vert z^{2} - w^{3} \vert^{-2\beta} . $$
The substitution $(z',w')=T(z,w)$ shows that
$$  I_{r+1} = 2^{(12 \beta- 10)} I_{r}. $$
Thus $\sum_{r} I_{r}$ is finite if $\beta<5/6$ and the union of the $\Omega_{r}$
cover $B^{4}\setminus \{0\}$.

\

Now we give the main proof. Let $x_{0}\in X$ be the point with stabiliser
$\Gamma$. We identify $SO(3)$-invariant functions on $X\setminus D$ with
$\Gamma$-invariant functions on $\cH=PSL(2,\bC)/SO(3)$ as in (4.2). The function $f_{0}$ on $x\setminus D$ corresponds to a convex
function $\phi_{0}$ on $\cH= PSL(2,\bC)/SO(3)$ which is an \lq\lq admissible
potential'' in the language of (4.2). For any other admissible potential $\phi$
the difference $\phi-\phi_{0}$ corresponds to $\psi$, restricted to $X\setminus
D$. The normalisation that $\max \psi=0$ becomes the condition that  $ {\rm
sup}\ \phi-\phi_{0}=0$,
and in particular $\phi\leq \phi_{0}$. 

Let $P_{0}\in \cH$ be the identity coset. It is the unique point fixed by
the action of $\Gamma$.
Any admissible potential function $\phi$ on $\cH$ is proper and bounded below so achieves
a minimum in $\cH$. By  the convexity and $\Gamma$-invariance this minimum
must occur at $P_{0}$.   Set $\phi(P_{0})=- b$. Then the inequality $\phi_{0}\geq
- b$ translates back into the statement that $\psi\geq f_{0} -b$. So
$$  \int_{Z} e^{-\beta \psi} d\mu_{0}\leq e^{b \beta} \int_{Z} f_{0}^{-\beta}d\mu_{0}.
$$
By Lemma 4, it suffices to obtain an upper bound on $b$. 
Let $B$ be the geodesic ball in $\cH$ centred on $P_{0}$, of radius $1$ say,
and let $ \oa$ be the maximum value of $\phi_{0}$ on $B$, so for any $\phi$
we have $\phi\leq  \oa$ on $B$. Convexity along geodesics emanating from
$P_{0}$ implies that
$$  \phi(Q)\leq -b+ (\oa+b)\ {\rm dist}(Q,P_{0}), $$
for any point $Q$ in $B$. In particular, on the ball $\frac{1}{2} B$ of radius $1/2$ about $P_{0}$ we have $\phi \leq (\oa-b)/2$. 

Take the inverse image in $PSL(2,\bC)$ of the ball $\frac{1}{2} B$ and map
this to $X$ by $g\rightarrow g(x_{0})$. The image obviously contains a neighbourhood
$N$ of $x_{0}$ and on $N$ we have $\psi\leq \frac{\oa-b}{2}+ f_{0}$. Then
Lemma 3 implies that
$b$ cannot be very large. In fact, if  the minimum of $f_{0}$ on $N$ is $
\ua$, we have $\psi\leq ( \frac{\oa}{2}- \ua)- \frac{b}{2}$ on $N$, so
$$ -M\leq \int_{N} \psi\  d\mu_{0} \leq ((\frac{\oa}{2}-\ua) - \frac{b}{2}) {\rm Vol}(N), $$
hence
$$ b\leq (\oa-2 \ua)+  \frac{2M}{{\rm Vol}(N)} $$
where $M$ is as in Lemma 3. This completes the proof of Theorem 3.

Notice that the same argument can be applied in the toric case, when the
polytope $P$ has a group $\Gamma$ of symmetries, as discussed in (4.2). We should
suppose that $\Gamma$ has a unique fixed point in $P$: then the proof proceeds
exactly as before.
 The analogue of Lemma 4 holds with $\beta< 1$
since the local models for the zeros of $s$ are $f_{p}(z_{1},
\dots, z_{n})=0$ where $f_{p}(z_{1},\dots, z_{n})=z_{1}\dots z_{p}$ and 
$\vert f_{p}\vert^{-2\beta}$ is locally integrable for $\beta<1$. the conclusion
is that
the $\alpha$-invariant in this case  is $1$.
which is a theorem of Batyrev and Selinova \cite{kn:BS}. (Song gave another
proof in \cite{kn:S}, and showed conversely that for polytopes which do not have such
a symmetry group the $\alpha$-invariant never exceeds $n/n+1$.)


\end{document}